\documentclass[a4paper,leqno]{article}
\usepackage{amsfonts}
\usepackage{amsmath}
\usepackage{amssymb}
\numberwithin{equation}{subsection}

\makeatletter
\def\diagram{\leftwidth=\z@ \rightwidth=\z@ \topheight=\z@
\botheight=\z@ \setbox\@picbox\hbox\bgroup}
\def\enddiagram{\egroup\wd\@picbox\rightwidth\unitlength
\ht\@picbox\topheight\unitlength \dp\@picbox\botheight\unitlength
\hskip\leftwidth\unitlength\box\@picbox}
\def\bfig{\begin{diagram}}
\def\efig{\end{diagram}}
\newcount\wideness \newcount\leftwidth \newcount\rightwidth
\newcount\highness \newcount\topheight \newcount\botheight
\def\ratchet#1#2{\ifnum#1<#2 \global #1=#2 \fi}
\def\putbox(#1,#2)#3{%
\horsize{\wideness}{#3} \divide\wideness by 2 {\advance\wideness
by #1 \ratchet{\rightwidth}{\wideness}} {\advance\wideness by -#1
\ratchet{\leftwidth}{\wideness}} \vertsize{\highness}{#3}
\divide\highness by 2 {\advance\highness by #2
\ratchet{\topheight}{\highness}} {\advance\highness by -#2
\ratchet{\botheight}{\highness}} \put(#1,#2){\makebox(0,0){$#3$}}}
\def\putlbox(#1,#2)#3{%
\horsize{\wideness}{#3} {\advance\wideness by #1
\ratchet{\rightwidth}{\wideness}} {\ratchet{\leftwidth}{-#1}}
\vertsize{\highness}{#3} \divide\highness by 2 {\advance\highness
by #2 \ratchet{\topheight}{\highness}} {\advance\highness by -#2
\ratchet{\botheight}{\highness}}
\put(#1,#2){\makebox(0,0)[l]{$#3$}}}
\def\putrbox(#1,#2)#3{%
\horsize{\wideness}{#3} {\ratchet{\rightwidth}{#1}}
{\advance\wideness by -#1 \ratchet{\leftwidth}{\wideness}}
\vertsize{\highness}{#3} \divide\highness by 2 {\advance\highness
by #2 \ratchet{\topheight}{\highness}} {\advance\highness by -#2
\ratchet{\botheight}{\highness}}
\put(#1,#2){\makebox(0,0)[r]{$#3$}}}

\def\adjust[#1]{} 
\newcount \coefa
\newcount \coefb
\newcount \coefc
\newcount\tempcounta
\newcount\tempcountb
\newcount\tempcountc
\newcount\tempcountd
\newcount\xext
\newcount\yext
\newcount\xoff
\newcount\yoff
\newcount\gap%
\newcount\arrowtypea
\newcount\arrowtypeb
\newcount\arrowtypec
\newcount\arrowtyped
\newcount\arrowtypee
\newcount\height
\newcount\width
\newcount\xpos
\newcount\ypos
\newcount\run
\newcount\rise
\newcount\arrowlength
\newcount\halflength
\newcount\arrowtype
\newdimen\tempdimen
\newdimen\xlen
\newdimen\ylen
\newsavebox{\tempboxa}%
\newsavebox{\tempboxb}%
\newsavebox{\tempboxc}%
\newdimen\w@dth
\def\setw@dth#1#2{\setbox\z@\hbox{$#1$}\w@dth=\wd\z@
\setbox\@ne\hbox{$#2$}\ifnum\w@dth<\wd\@ne \w@dth=\wd\@ne \fi
\advance\w@dth by 1.2em}
\def\t@^#1_#2{\def\n@one{#1}\def\n@two{#2}\mathrel{\setw@dth{#1}{#2}
\mathop{\hbox to \w@dth{\rightarrowfill}}\limits
\ifx\n@one\empty\else ^{\box\z@}\fi \ifx\n@two\empty\else
_{\box\@ne}\fi}}
\def\t@@^#1{\@ifnextchar_ {\t@^{#1}}{\t@^{#1}_{}}}
\def\to{\@ifnextchar^ {\t@@}{\t@@^{}}}
\def\t@left^#1_#2{\def\n@one{#1}\def\n@two{#2}\mathrel{\setw@dth{#1}{#2}
\mathop{\hbox to \w@dth{\leftarrowfill}}\limits
\ifx\n@one\empty\else ^{\box\z@}\fi \ifx\n@two\empty\else
_{\box\@ne}\fi}}
\def\t@@left^#1{\@ifnextchar_ {\t@left^{#1}}{\t@left^{#1}_{}}}
\def\toleft{\@ifnextchar^ {\t@@left}{\t@@left^{}}}
\def\two@^#1_#2{\def\n@one{#1}\def\n@two{#2}\mathrel{\setw@dth{#1}{#2}
\mathop{\vcenter{\hbox to \w@dth{\rightarrowfill}\kern-1.7ex
                 \hbox to \w@dth{\rightarrowfill}}%
       }\limits
\ifx\n@one\empty\else ^{\box\z@}\fi \ifx\n@two\empty\else
_{\box\@ne}\fi}}
\def\tw@@^#1{\@ifnextchar_ {\two@^{#1}}{\two@^{#1}_{}}}
\def\two{\@ifnextchar^ {\tw@@}{\tw@@^{}}}
\def\tofr@^#1_#2{\def\n@one{#1}\def\n@two{#2}\mathrel{\setw@dth{#1}{#2}
\mathop{\vcenter{\hbox to \w@dth{\rightarrowfill}\kern-1.7ex
                 \hbox to \w@dth{\leftarrowfill}}%
       }\limits
\ifx\n@one\empty\else ^{\box\z@}\fi \ifx\n@two\empty\else
_{\box\@ne}\fi}}
\def\t@fr@^#1{\@ifnextchar_ {\tofr@^{#1}}{\tofr@^{#1}_{}}}
\def\tofro{\@ifnextchar^ {\t@fr@}{\t@fr@^{}}}

\def\mon{\mathop{\m@th\hbox to
      14.6\P@{\lasyb\char'51\hskip-2.1\P@$\arrext$\hss
$\mathord\rightarrow$}}\limits} 
\def\leftmono{\mathrel{\m@th\hbox to
14.6\P@{$\mathord\leftarrow$\hss$\arrext$\hskip-2.1\P@\lasyb\char'50%
}}\limits} 
\mathchardef\arrext="0200       

\def\settypes(#1,#2,#3){\arrowtypea#1 \arrowtypeb#2 \arrowtypec#3}
\def\settoheight#1#2{\setbox\@tempboxa\hbox{#2}#1\ht\@tempboxa\relax}%
\def\settodepth#1#2{\setbox\@tempboxa\hbox{#2}#1\dp\@tempboxa\relax}%
\def\settokens[#1`#2`#3`#4]{%
     \def\tokena{#1}\def\tokenb{#2}\def\tokenc{#3}\def\tokend{#4}}
\def\setsqparms[#1`#2`#3`#4;#5`#6]{%
\arrowtypea #1 \arrowtypeb #2 \arrowtypec #3 \arrowtyped #4
\width #5 \height #6 }
\def\setpos(#1,#2){\xpos=#1 \ypos#2}

\def\settriparms[#1`#2`#3;#4]{\settripairparms[#1`#2`#3`1`1;#4]}%
\def\settripairparms[#1`#2`#3`#4`#5;#6]{%
\arrowtypea #1 \arrowtypeb #2 \arrowtypec #3 \arrowtyped #4
\arrowtypee #5 \width #6 \height #6 }
\def\resetparms{\settripairparms[1`1`1`1`1;500]\width 500}
\resetparms
\def\mvector(#1,#2)#3{
\put(0,0){\vector(#1,#2){#3}}%
\put(0,0){\vector(#1,#2){26}}%
}
\def\evector(#1,#2)#3{{
\arrowlength #3
\put(0,0){\vector(#1,#2){\arrowlength}}%
\advance \arrowlength by-30
\put(0,0){\vector(#1,#2){\arrowlength}}%
}}
\def\horsize#1#2{%
\settowidth{\tempdimen}{$#2$}%
#1=\tempdimen \divide #1 by\unitlength }
\def\vertsize#1#2{%
\settoheight{\tempdimen}{$#2$}%
#1=\tempdimen
\settodepth{\tempdimen}{$#2$}%
\advance #1 by\tempdimen \divide #1 by\unitlength }
\def\putvector(#1,#2)(#3,#4)#5#6{{%
\ifnum3<\arrowtype \putdashvector(#1,#2)(#3,#4)#5\arrowtype \else
\ifnum\arrowtype<-3 \putdashvector(#1,#2)(#3,#4)#5\arrowtype \else
\xpos=#1 \ypos=#2 \run=#3 \rise=#4 \arrowlength=#5 \ifnum
\arrowtype<0
    \ifnum \run=0
        \advance \ypos by-\arrowlength
    \else
        \tempcounta \arrowlength
        \multiply \tempcounta by\rise
        \divide \tempcounta by\run
        \ifnum\run>0
            \advance \xpos by\arrowlength
            \advance \ypos by\tempcounta
        \else
            \advance \xpos by-\arrowlength
            \advance \ypos by-\tempcounta
        \fi
    \fi
    \multiply \arrowtype by-1
    \multiply \rise by-1
    \multiply \run by-1
\fi \ifcase \arrowtype
\or \put(\xpos,\ypos){\vector(\run,\rise){\arrowlength}}%
\or \put(\xpos,\ypos){\mvector(\run,\rise)\arrowlength}%
\or \put(\xpos,\ypos){\evector(\run,\rise){\arrowlength}}%
\fi\fi\fi }}
\def\putsplitvector(#1,#2)#3#4{
\xpos #1 \ypos #2 \arrowtype #4 \halflength #3 \arrowlength #3
\gap 140 \advance \halflength by-\gap \divide \halflength by2
\ifnum\arrowtype>0
   \ifcase \arrowtype
   \or \put(\xpos,\ypos){\line(0,-1){\halflength}}%
       \advance\ypos by-\halflength
       \advance\ypos by-\gap
       \put(\xpos,\ypos){\vector(0,-1){\halflength}}%
   \or \put(\xpos,\ypos){\line(0,-1)\halflength}%
       \put(\xpos,\ypos){\vector(0,-1)3}%
       \advance\ypos by-\halflength
       \advance\ypos by-\gap
       \put(\xpos,\ypos){\vector(0,-1){\halflength}}%
   \or \put(\xpos,\ypos){\line(0,-1)\halflength}%
       \advance\ypos by-\halflength
       \advance\ypos by-\gap
       \put(\xpos,\ypos){\evector(0,-1){\halflength}}%
   \fi
\else \arrowtype=-\arrowtype
   \ifcase\arrowtype
   \or \advance \ypos by-\arrowlength
       \put(\xpos,\ypos){\line(0,1){\halflength}}%
       \advance\ypos by\halflength
       \advance\ypos by\gap
       \put(\xpos,\ypos){\vector(0,1){\halflength}}%
   \or \advance \ypos by-\arrowlength
       \put(\xpos,\ypos){\line(0,1)\halflength}%
       \put(\xpos,\ypos){\vector(0,1)3}%
       \advance\ypos by\halflength
       \advance\ypos by\gap
       \put(\xpos,\ypos){\vector(0,1){\halflength}}%
   \or \advance \ypos by-\arrowlength
       \put(\xpos,\ypos){\line(0,1)\halflength}%
       \advance\ypos by\halflength
       \advance\ypos by\gap
       \put(\xpos,\ypos){\evector(0,1){\halflength}}%
   \fi
\fi }
\def\putmorphism(#1)(#2,#3)[#4`#5`#6]#7#8#9{{%
\run #2 \rise #3 \ifnum\rise=0
  \puthmorphism(#1)[#4`#5`#6]{#7}{#8}#9%
\else\ifnum\run=0
  \putvmorphism(#1)[#4`#5`#6]{#7}{#8}#9%
\else
\setpos(#1)%
\arrowlength #7 \arrowtype #8 \ifnum\run=0 \else\ifnum\rise=0
\else \ifnum\run>0
    \coefa=1
\else
   \coefa=-1
\fi \ifnum\arrowtype>0
   \coefb=0
   \coefc=-1
\else
   \coefb=\coefa
   \coefc=1
   \arrowtype=-\arrowtype
\fi \width=2 \multiply \width by\run \divide \width by\rise
\ifnum \width<0  \width=-\width\fi \advance\width by60 \if l#9
\width=-\width\fi
\putbox(\xpos,\ypos){#4}
{\multiply \coefa by\arrowlength
\advance\xpos by\coefa \multiply \coefa by\rise \divide \coefa
by\run \advance \ypos by\coefa
\putbox(\xpos,\ypos){#5} }%
{\multiply \coefa by\arrowlength
\divide \coefa by2 \advance \xpos by\coefa \advance \xpos by\width
\multiply \coefa by\rise \divide \coefa by\run \advance \ypos
by\coefa
\if l#9%
   \putrbox(\xpos,\ypos){#6}%
\else\if r#9%
   \putlbox(\xpos,\ypos){#6}%
\fi\fi }%
{\multiply \rise by-\coefc
\multiply \run by-\coefc \multiply \coefb by\arrowlength \advance
\xpos by\coefb \multiply \coefb by\rise \divide \coefb by\run
\advance \ypos by\coefb \multiply \coefc by70 \advance \ypos
by\coefc \multiply \coefc by\run \divide \coefc by\rise \advance
\xpos by\coefc \multiply \coefa by140 \multiply \coefa by\run
\divide \coefa by\rise \advance \arrowlength by\coefa
\ifcase\arrowtype
\or \put(\xpos,\ypos){\vector(\run,\rise){\arrowlength}}%
\or \put(\xpos,\ypos){\mvector(\run,\rise){\arrowlength}}%
\or \put(\xpos,\ypos){\evector(\run,\rise){\arrowlength}}%
\fi}\fi\fi\fi\fi}}

\newcount\numbdashes \newcount\lengthdash \newcount\increment
\def\howmanydashes{
\numbdashes=\arrowlength \lengthdash=40 \divide\numbdashes by
\lengthdash \lengthdash=\arrowlength \divide\lengthdash by
\numbdashes
\increment=\lengthdash \multiply\lengthdash by 3
\divide\lengthdash by 5 }
\def\putdashvector(#1)(#2,#3)#4#5{%
\ifnum#3=0 \putdashhvector(#1){#4}#5 \else \ifnum#2=0
\putdashvvector(#1){#4}#5\fi\fi}
\def\putdashhvector(#1,#2)#3#4{{%
\arrowlength=#3 \howmanydashes
\multiput(#1,#2)(\increment,0){\numbdashes}%
{\vrule height .4pt width \lengthdash\unitlength} \arrowtype=#4
\xpos=#1 \ifnum\arrowtype<0 \advance\arrowtype by 7 \fi
\ifcase\arrowtype \or \advance\xpos by 10
    \put(\xpos,#2){\vector(-1,0){\lengthdash}}
    \advance\xpos by 40
    \put(\xpos,#2){\vector(-1,0){\lengthdash}}
\or \advance \xpos by 10
    \put(\xpos,#2){\vector(-1,0){\lengthdash}}
    \advance\xpos by  \arrowlength
    \advance\xpos by  -50
    \put(\xpos,#2){\vector(-1,0){\lengthdash}}
\or \advance\xpos by 10
    \put(\xpos,#2){\vector(-1,0){\lengthdash}}
\or \advance\xpos by \arrowlength
    \advance\xpos by -\lengthdash
    \put(\xpos,#2){\vector(1,0){\lengthdash}}
\or {\advance\xpos by 10
    \put(\xpos,#2){\vector(1,0){\lengthdash}}}
    \advance\xpos by \arrowlength
    \advance\xpos by -\lengthdash
    \put(\xpos,#2){\vector(1,0){\lengthdash}}
\or \advance\xpos by \arrowlength
    \advance\xpos by -\lengthdash
    \put(\xpos,#2){\vector(1,0){\lengthdash}}
    \advance\xpos by -40
    \put(\xpos,#2){\vector(1,0){\lengthdash}}
   \fi
}}
\def\putdashvvector(#1,#2)#3#4{{%
\arrowlength=#3 \howmanydashes \ypos=#2 \advance\ypos by
-\arrowlength
\multiput(#1,#2)(0,\increment){\numbdashes}%
    {\vrule width .4pt height \lengthdash\unitlength}
\arrowtype=#4 \ypos=#2 \ifnum\arrowtype<0 \advance\arrowtype by 7
\fi \ifcase\arrowtype \or \advance\ypos by \arrowlength
\advance\ypos by -40
    \put(#1,\ypos){\vector(0,1){\lengthdash}}
    \advance\ypos by -40
    \put(#1,\ypos){\vector(0,1){\lengthdash}}
\or \advance\ypos by 10
    \put(#1,\ypos){\vector(0,1){\lengthdash}}
    \advance\ypos by \arrowlength \advance\ypos by -40
    \put(#1,\ypos){\vector(0,1){\lengthdash}}
\or \advance\ypos by \arrowlength \advance\ypos by -40
    \put(#1,\ypos){\vector(0,1){\lengthdash}}
\or \advance\ypos by 10
    \put(#1,\ypos){\vector(0,-1){\lengthdash}}
\or \advance\ypos by 10
    \put(#1,\ypos){\vector(0,-1){\lengthdash}}
    \advance\ypos by \arrowlength \advance\ypos by -40
    \put(#1,\ypos){\vector(0,-1){\lengthdash}}
\or \advance\ypos by 10
    \put(#1,\ypos){\vector(0,-1){\lengthdash}}
    \advance\ypos by 40
    \put(#1,\ypos){\vector(0,-1){\lengthdash}}
\fi }}
\def\puthmorphism(#1,#2)[#3`#4`#5]#6#7#8{{%
\xpos #1 \ypos #2 \width #6 \arrowlength #6 \arrowtype=#7
\putbox(\xpos,\ypos){#3\vphantom{#4}}%
{\advance \xpos by\arrowlength
\putbox(\xpos,\ypos){\vphantom{#3}#4}}%
\horsize{\tempcounta}{#3}%
\horsize{\tempcountb}{#4}%
\divide \tempcounta by2 \divide \tempcountb by2 \advance
\tempcounta by30 \advance \tempcountb by30 \advance \xpos
by\tempcounta \advance \arrowlength by-\tempcounta \advance
\arrowlength by-\tempcountb
\putvector(\xpos,\ypos)(1,0)\arrowlength\arrowtype \divide
\arrowlength by2 \advance \xpos by\arrowlength
\vertsize{\tempcounta}{#5}%
\divide\tempcounta by2 \advance \tempcounta by20
\if a#8 %
   \advance \ypos by\tempcounta
   \putbox(\xpos,\ypos){#5}%
\else
   \advance \ypos by-\tempcounta
   \putbox(\xpos,\ypos){#5}%
\fi}}
\def\putvmorphism(#1,#2)[#3`#4`#5]#6#7#8{{%
\xpos #1 \ypos #2 \arrowlength #6 \arrowtype #7
\settowidth{\xlen}{$#5$}%
\putbox(\xpos,\ypos){#3}%
{\advance \ypos by-\arrowlength
\putbox(\xpos,\ypos){#4}}%
{\advance\arrowlength by-140 \advance \ypos by-70 \ifdim\xlen>0pt
   \if m#8%
      \putsplitvector(\xpos,\ypos)\arrowlength\arrowtype
   \else
   \putvector(\xpos,\ypos)(0,-1)\arrowlength\arrowtype
   \fi
\else
   \putvector(\xpos,\ypos)(0,-1)\arrowlength\arrowtype
\fi}%
\ifdim\xlen>0pt
   \divide \arrowlength by2
   \advance\ypos by-\arrowlength
   \if l#8%
      \advance \xpos by-40
      \putrbox(\xpos,\ypos){#5}%
   \else\if r#8%
      \advance \xpos by40
      \putlbox(\xpos,\ypos){#5}%
   \else
      \putbox(\xpos,\ypos){#5}%
   \fi\fi
\fi }}
\def\putsquarep<#1>(#2)[#3;#4`#5`#6`#7]{{%
\setsqparms[#1]%
\setpos(#2)%
\settokens[#3]%
\puthmorphism(\xpos,\ypos)[\tokenc`\tokend`{#7}]{\width}{\arrowtyped}b%
\advance\ypos by \height
\puthmorphism(\xpos,\ypos)[\tokena`\tokenb`{#4}]{\width}{\arrowtypea}a%
\putvmorphism(\xpos,\ypos)[``{#5}]{\height}{\arrowtypeb}l%
\advance\xpos by \width
\putvmorphism(\xpos,\ypos)[``{#6}]{\height}{\arrowtypec}r%
}}
\def\putsquare{\@ifnextchar <{\putsquarep}{\putsquarep%
   <\arrowtypea`\arrowtypeb`\arrowtypec`\arrowtyped;\width`\height>}}
\def\square{\@ifnextchar< {\squarep}{\squarep
   <\arrowtypea`\arrowtypeb`\arrowtypec`\arrowtyped;\width`\height>}}
\def\squarep<#1>[#2`#3`#4`#5;#6`#7`#8`#9]{{
\setsqparms[#1]
\diagram
\putsquarep<\arrowtypea`\arrowtypeb`\arrowtypec`
\arrowtyped;\width`\height>
(0,0)[#2`#3`#4`{#5};#6`#7`#8`{#9}]
\enddiagram
}}                                                 
\def\putptrianglep<#1>(#2,#3)[#4`#5`#6;#7`#8`#9]{{%
\settriparms[#1]%
\xpos=#2 \ypos=#3 \advance\ypos by \height
\puthmorphism(\xpos,\ypos)[#4`#5`{#7}]{\height}{\arrowtypea}a%
\putvmorphism(\xpos,\ypos)[`#6`{#8}]{\height}{\arrowtypeb}l%
\advance\xpos by\height
\putmorphism(\xpos,\ypos)(-1,-1)[``{#9}]{\height}{\arrowtypec}r%
}}
\def\putptriangle{\@ifnextchar <{\putptrianglep}{\putptrianglep
   <\arrowtypea`\arrowtypeb`\arrowtypec;\height>}}
\def\ptriangle{\@ifnextchar <{\ptrianglep}{\ptrianglep
   <\arrowtypea`\arrowtypeb`\arrowtypec;\height>}}
\def\ptrianglep<#1>[#2`#3`#4;#5`#6`#7]{{
\settriparms[#1]
\diagram
\putptrianglep<\arrowtypea`\arrowtypeb`
\arrowtypec;\height>
(0,0)[#2`#3`#4;#5`#6`{#7}]
\enddiagram
}}                                            
\def\putqtrianglep<#1>(#2,#3)[#4`#5`#6;#7`#8`#9]{{%
\settriparms[#1]%
\xpos=#2 \ypos=#3 \advance\ypos by\height
\puthmorphism(\xpos,\ypos)[#4`#5`{#7}]{\height}{\arrowtypea}a%
\putmorphism(\xpos,\ypos)(1,-1)[``{#8}]{\height}{\arrowtypeb}l%
\advance\xpos by\height
\putvmorphism(\xpos,\ypos)[`#6`{#9}]{\height}{\arrowtypec}r%
}}
\def\putqtriangle{\@ifnextchar <{\putqtrianglep}{\putqtrianglep
   <\arrowtypea`\arrowtypeb`\arrowtypec;\height>}}
\def\qtriangle{\@ifnextchar <{\qtrianglep}{\qtrianglep
   <\arrowtypea`\arrowtypeb`\arrowtypec;\height>}}
\def\qtrianglep<#1>[#2`#3`#4;#5`#6`#7]{{
\settriparms[#1]
\width=\height                                
\diagram
\putqtrianglep<\arrowtypea`\arrowtypeb`
\arrowtypec;\height>
(0,0)[#2`#3`#4;#5`#6`{#7}]
\enddiagram
}}
\def\putdtrianglep<#1>(#2,#3)[#4`#5`#6;#7`#8`#9]{{%
\settriparms[#1]%
\xpos=#2 \ypos=#3
\puthmorphism(\xpos,\ypos)[#5`#6`{#9}]{\height}{\arrowtypec}b%
\advance\xpos by \height \advance\ypos by\height
\putmorphism(\xpos,\ypos)(-1,-1)[``{#7}]{\height}{\arrowtypea}l%
\putvmorphism(\xpos,\ypos)[#4``{#8}]{\height}{\arrowtypeb}r%
}}
\def\putdtriangle{\@ifnextchar <{\putdtrianglep}{\putdtrianglep
   <\arrowtypea`\arrowtypeb`\arrowtypec;\height>}}
\def\dtriangle{\@ifnextchar <{\dtrianglep}{\dtrianglep
   <\arrowtypea`\arrowtypeb`\arrowtypec;\height>}}
\def\dtrianglep<#1>[#2`#3`#4;#5`#6`#7]{{
\settriparms[#1]
\width=\height                                
\diagram
\putdtrianglep<\arrowtypea`\arrowtypeb`
\arrowtypec;\height>
(0,0)[#2`#3`#4;#5`#6`{#7}]
\enddiagram
}}
\def\putbtrianglep<#1>(#2,#3)[#4`#5`#6;#7`#8`#9]{{%
\settriparms[#1]%
\xpos=#2 \ypos=#3
\puthmorphism(\xpos,\ypos)[#5`#6`{#9}]{\height}{\arrowtypec}b%
\advance\ypos by\height
\putmorphism(\xpos,\ypos)(1,-1)[``{#8}]{\height}{\arrowtypeb}r%
\putvmorphism(\xpos,\ypos)[#4``{#7}]{\height}{\arrowtypea}l%
}}
\def\putbtriangle{\@ifnextchar <{\putbtrianglep}{\putbtrianglep
   <\arrowtypea`\arrowtypeb`\arrowtypec;\height>}}
\def\btriangle{\@ifnextchar <{\btrianglep}{\btrianglep
   <\arrowtypea`\arrowtypeb`\arrowtypec;\height>}}
\def\btrianglep<#1>[#2`#3`#4;#5`#6`#7]{{
\settriparms[#1]
\width=\height                               
\diagram
\putbtrianglep<\arrowtypea`\arrowtypeb`
\arrowtypec;\height>
(0,0)[#2`#3`#4;#5`#6`{#7}]
\enddiagram
}}
\def\putAtrianglep<#1>(#2,#3)[#4`#5`#6;#7`#8`#9]{{%
\settriparms[#1]%
\xpos=#2 \ypos=#3 {\multiply \height by2
\puthmorphism(\xpos,\ypos)[#5`#6`{#9}]{\height}{\arrowtypec}b}%
\advance\xpos by\height \advance\ypos by\height
\putmorphism(\xpos,\ypos)(-1,-1)[#4``{#7}]{\height}{\arrowtypea}l%
\putmorphism(\xpos,\ypos)(1,-1)[``{#8}]{\height}{\arrowtypeb}r%
}}
\def\putAtriangle{\@ifnextchar <{\putAtrianglep}{\putAtrianglep
   <\arrowtypea`\arrowtypeb`\arrowtypec;\height>}}
\def\Atriangle{\@ifnextchar <{\Atrianglep}{\Atrianglep
   <\arrowtypea`\arrowtypeb`\arrowtypec;\height>}}
\def\Atrianglep<#1>[#2`#3`#4;#5`#6`#7]{{
\settriparms[#1]
\width=\height                                     
\diagram
\putAtrianglep<\arrowtypea`\arrowtypeb`
\arrowtypec;\height>
(0,0)[#2`#3`#4;#5`#6`{#7}]
\enddiagram
}}
\def\putAtrianglepairp<#1>(#2)[#3;#4`#5`#6`#7`#8]{{%
\settripairparms[#1]%
\setpos(#2)%
\settokens[#3]%
\puthmorphism(\xpos,\ypos)[\tokenb`\tokenc`{#7}]{\height}{\arrowtyped}b%
\advance\xpos by\height
\puthmorphism(\xpos,\ypos)[\phantom{\tokenc}`\tokend`{#8}]%
{\height}{\arrowtypee}b%
\advance\ypos by\height
\putmorphism(\xpos,\ypos)(-1,-1)[\tokena``{#4}]{\height}{\arrowtypea}l%
\putvmorphism(\xpos,\ypos)[``{#5}]{\height}{\arrowtypeb}m%
\putmorphism(\xpos,\ypos)(1,-1)[``{#6}]{\height}{\arrowtypec}r%
}}
\def\putAtrianglepair{\@ifnextchar <{\putAtrianglepairp}{\putAtrianglepairp%
   <\arrowtypea`\arrowtypeb`\arrowtypec`\arrowtyped`\arrowtypee;\height>}}
\def\Atrianglepair{\@ifnextchar <{\Atrianglepairp}{\Atrianglepairp%
   <\arrowtypea`\arrowtypeb`\arrowtypec`\arrowtyped`\arrowtypee;\height>}}
\def\Atrianglepairp<#1>[#2;#3`#4`#5`#6`#7]{{
\settripairparms[#1]
\settokens[#2]
\width=\height                                
\diagram
\putAtrianglepairp                            
<\arrowtypea`\arrowtypeb`\arrowtypec`
\arrowtyped`\arrowtypee;\height>
(0,0)[{#2};#3`#4`#5`#6`{#7}]
\enddiagram
}}
\def\putVtrianglep<#1>(#2,#3)[#4`#5`#6;#7`#8`#9]{{%
\settriparms[#1]%
\xpos=#2 \ypos=#3 \advance\ypos by\height {\multiply\height by2
\puthmorphism(\xpos,\ypos)[#4`#5`{#7}]{\height}{\arrowtypea}a}%
\putmorphism(\xpos,\ypos)(1,-1)[`#6`{#8}]{\height}{\arrowtypeb}l%
\advance\xpos by\height \advance\xpos by\height
\putmorphism(\xpos,\ypos)(-1,-1)[``{#9}]{\height}{\arrowtypec}r%
}}
\def\putVtriangle{\@ifnextchar <{\putVtrianglep}{\putVtrianglep
   <\arrowtypea`\arrowtypeb`\arrowtypec;\height>}}
\def\Vtriangle{\@ifnextchar <{\Vtrianglep}{\Vtrianglep
   <\arrowtypea`\arrowtypeb`\arrowtypec;\height>}}
\def\Vtrianglep<#1>[#2`#3`#4;#5`#6`#7]{{
\settriparms[#1]
\width=\height                                 
\diagram
\putVtrianglep<\arrowtypea`\arrowtypeb`
\arrowtypec;\height>
(0,0)[#2`#3`#4;#5`#6`{#7}]
\enddiagram
}}
\def\putVtrianglepairp<#1>(#2)[#3;#4`#5`#6`#7`#8]{{
\settripairparms[#1]%
\setpos(#2)%
\settokens[#3]%
\advance\ypos by\height
\putmorphism(\xpos,\ypos)(1,-1)[`\tokend`{#6}]{\height}{\arrowtypec}l%
\puthmorphism(\xpos,\ypos)[\tokena`\tokenb`{#4}]{\height}{\arrowtypea}a%
\advance\xpos by\height
\puthmorphism(\xpos,\ypos)[\phantom{\tokenb}`\tokenc`{#5}]%
{\height}{\arrowtypeb}a%
\putvmorphism(\xpos,\ypos)[``{#7}]{\height}{\arrowtyped}m%
\advance\xpos by\height
\putmorphism(\xpos,\ypos)(-1,-1)[``{#8}]{\height}{\arrowtypee}r%
}}
\def\putVtrianglepair{\@ifnextchar <{\putVtrianglepairp}{\putVtrianglepairp%
    <\arrowtypea`\arrowtypeb`\arrowtypec`\arrowtyped`\arrowtypee;\height>}}
\def\Vtrianglepair{\@ifnextchar <{\Vtrianglepairp}{\Vtrianglepairp%
    <\arrowtypea`\arrowtypeb`\arrowtypec`\arrowtyped`\arrowtypee;\height>}}
\def\Vtrianglepairp<#1>[#2;#3`#4`#5`#6`#7]{{
\settripairparms[#1]
\settokens[#2]
\diagram
\putVtrianglepairp                             
<\arrowtypea`\arrowtypeb`\arrowtypec`
\arrowtyped`\arrowtypee;\height>
(0,0)[{#2};#3`#4`#5`#6`{#7}]
\enddiagram
}}

\def\putCtrianglep<#1>(#2,#3)[#4`#5`#6;#7`#8`#9]{{%
\settriparms[#1]%
\xpos=#2 \ypos=#3 \advance\ypos by\height
\putmorphism(\xpos,\ypos)(1,-1)[``{#9}]{\height}{\arrowtypec}l%
\advance\xpos by\height \advance\ypos by\height
\putmorphism(\xpos,\ypos)(-1,-1)[#4`#5`{#7}]{\height}{\arrowtypea}l%
{\multiply\height by 2
\putvmorphism(\xpos,\ypos)[`#6`{#8}]{\height}{\arrowtypeb}r}%
}}
\def\putCtriangle{\@ifnextchar <{\putCtrianglep}{\putCtrianglep
    <\arrowtypea`\arrowtypeb`\arrowtypec;\height>}}
\def\Ctriangle{\@ifnextchar <{\Ctrianglep}{\Ctrianglep
    <\arrowtypea`\arrowtypeb`\arrowtypec;\height>}}
\def\Ctrianglep<#1>[#2`#3`#4;#5`#6`#7]{{
\settriparms[#1]
\width=\height                               
\diagram
\putCtrianglep<\arrowtypea`\arrowtypeb`
\arrowtypec;\height>
(0,0)[#2`#3`#4;#5`#6`{#7}]
\enddiagram
}}                                           
\def\putDtrianglep<#1>(#2,#3)[#4`#5`#6;#7`#8`#9]{{%
\settriparms[#1]%
\xpos=#2 \ypos=#3 \advance\xpos by\height \advance\ypos by\height
\putmorphism(\xpos,\ypos)(-1,-1)[``{#9}]{\height}{\arrowtypec}r%
\advance\xpos by-\height \advance\ypos by\height
\putmorphism(\xpos,\ypos)(1,-1)[`#5`{#8}]{\height}{\arrowtypeb}r%
{\multiply\height by 2
\putvmorphism(\xpos,\ypos)[#4`#6`{#7}]{\height}{\arrowtypea}l}%
}}
\def\putDtriangle{\@ifnextchar <{\putDtrianglep}{\putDtrianglep
    <\arrowtypea`\arrowtypeb`\arrowtypec;\height>}}
\def\Dtriangle{\@ifnextchar <{\Dtrianglep}{\Dtrianglep
   <\arrowtypea`\arrowtypeb`\arrowtypec;\height>}}
\def\Dtrianglep<#1>[#2`#3`#4;#5`#6`#7]{{
\settriparms[#1]
\width=\height                              
\diagram
\putDtrianglep<\arrowtypea`\arrowtypeb`
\arrowtypec;\height>
(0,0)[#2`#3`#4;#5`#6`{#7}]
\enddiagram
}}                                          
\def\setrecparms[#1`#2]{\width=#1 \height=#2}%
\def\recursep<#1`#2>[#3;#4`#5`#6`#7`#8]{{%
\width=#1 \height=#2 \settokens[#3]
\settowidth{\tempdimen}{$\tokena$} \ifdim\tempdimen=0pt
  \savebox{\tempboxa}{\hbox{$\tokenb$}}%
  \savebox{\tempboxb}{\hbox{$\tokend$}}%
  \savebox{\tempboxc}{\hbox{$#6$}}%
\else
  \savebox{\tempboxa}{\hbox{$\hbox{$\tokena$}\times\hbox{$\tokenb$}$}}%
  \savebox{\tempboxb}{\hbox{$\hbox{$\tokena$}\times\hbox{$\tokend$}$}}%
  \savebox{\tempboxc}{\hbox{$\hbox{$\tokena$}\times\hbox{$#6$}$}}%
\fi \ypos=\height \divide\ypos by 2 \xpos=\ypos \advance\xpos by
\width \bfig
\putCtrianglep<-1`1`1;\ypos>(0,0)[`\tokenc`;#5`#6`{#7}]%
\puthmorphism(\ypos,0)[\tokend`\usebox{\tempboxb}`{#8}]{\width}{-1}b%
\puthmorphism(\ypos,\height)[\tokenb`\usebox{\tempboxa}`{#4}]{\width}{-1}a%
\advance\ypos by \width
\putvmorphism(\ypos,\height)[``\usebox{\tempboxc}]{\height}1r%
\efig }}
\def\recurse{\@ifnextchar <{\recursep}{\recursep<\width`\height>}}
\def\puttwohmorphisms(#1,#2)[#3`#4;#5`#6]#7#8#9{{%
\puthmorphism(#1,#2)[#3`#4`]{#7}0a \ypos=#2 \advance\ypos by 20
\puthmorphism(#1,\ypos)[\phantom{#3}`\phantom{#4}`#5]{#7}{#8}a
\advance\ypos by -40
\puthmorphism(#1,\ypos)[\phantom{#3}`\phantom{#4}`#6]{#7}{#9}b }}
\def\puttwovmorphisms(#1,#2)[#3`#4;#5`#6]#7#8#9{{%
\putvmorphism(#1,#2)[#3`#4`]{#7}0a \xpos=#1 \advance\xpos by -20
\putvmorphism(\xpos,#2)[\phantom{#3}`\phantom{#4}`#5]{#7}{#8}l
\advance\xpos by 40
\putvmorphism(\xpos,#2)[\phantom{#3}`\phantom{#4}`#6]{#7}{#9}r }}
\def\puthcoequalizer(#1)[#2`#3`#4;#5`#6`#7]#8#9{{%
\setpos(#1)%
\puttwohmorphisms(\xpos,\ypos)[#2`#3;#5`#6]{#8}11%
\advance\xpos by #8
\puthmorphism(\xpos,\ypos)[\phantom{#3}`#4`#7]{#8}1{#9} }}
\def\putvcoequalizer(#1)[#2`#3`#4;#5`#6`#7]#8#9{{%
\setpos(#1)%
\puttwovmorphisms(\xpos,\ypos)[#2`#3;#5`#6]{#8}11%
\advance\ypos by -#8
\putvmorphism(\xpos,\ypos)[\phantom{#3}`#4`#7]{#8}1{#9} }}
\def\putthreehmorphisms(#1)[#2`#3;#4`#5`#6]#7(#8)#9{{%
\setpos(#1) \settypes(#8)
\if a#9 %
     \vertsize{\tempcounta}{#5}%
     \vertsize{\tempcountb}{#6}%
     \ifnum \tempcounta<\tempcountb \tempcounta=\tempcountb \fi
\else
     \vertsize{\tempcounta}{#4}%
     \vertsize{\tempcountb}{#5}%
     \ifnum \tempcounta<\tempcountb \tempcounta=\tempcountb \fi
\fi \advance \tempcounta by 60
\puthmorphism(\xpos,\ypos)[#2`#3`#5]{#7}{\arrowtypeb}{#9}
\advance\ypos by \tempcounta
\puthmorphism(\xpos,\ypos)[\phantom{#2}`\phantom{#3}`#4]{#7}{\arrowtypea}{#9}
\advance\ypos by -\tempcounta \advance\ypos by -\tempcounta
\puthmorphism(\xpos,\ypos)[\phantom{#2}`\phantom{#3}`#6]{#7}{\arrowtypec}{#9}
}}
\def\setarrowtoks[#1`#2`#3`#4`#5`#6]{%
\def\toka{#1}
\def\tokb{#2}
\def\tokc{#3}
\def\tokd{#4}
\def\toke{#5}
\def\tokf{#6}
}
\def\hex{\@ifnextchar <{\hexp}{\hexp<1000`400>}}
\def\hexp<#1`#2>[#3`#4`#5`#6`#7`#8;#9]{%
\setarrowtoks[#9] \yext=#2 \advance \yext by #2 \xext=#1
\advance\xext by \yext \bfig
\putCtriangle<-1`0`1;#2>(0,0)[`#5`;\tokb``\tokd] \xext=#1
\yext=#2 \advance \yext by #2
\putsquare<1`0`0`1;\xext`\yext>(#2,0)[#3`#4`#7`#8;\toka```\tokf]
\advance \xext by #2
\putDtriangle<0`1`-1;#2>(\xext,0)[`#6`;`\tokc`\toke] \efig }

\begin{document}
\newtheorem{theorem}{Theorem}[section]
\newtheorem{lemma}[theorem]{Lemma}
\newtheorem{corollary}[theorem]{Corollary}
\newtheorem{conjecture}[theorem]{Conjecture}
\newtheorem{remark}[theorem]{Remark}
\newtheorem{condition}[theorem]{Condition}{\it}{\rm}
\newtheorem{definition}[theorem]{Definition}
\newtheorem{problem}[theorem]{Problem}
\newtheorem{example}[theorem]{Example}
\newtheorem{proposition}[theorem]{Proposition}
\newcommand{\cA}{{\mathcal A} }
\newcommand{\cB}{{\mathcal B} }
\newcommand{\cD}{{\mathcal D} }
\newcommand{\cE}{{\mathcal E} }
\newcommand{\cF}{{\mathcal F} }
\newcommand{\cG}{{\mathcal G} }
\newcommand{\cH}{{\mathcal H} }
\newcommand{\cI}{{\mathcal I} }
\newcommand{\cJ}{{\mathcal J} }
\newcommand{\cK}{{\mathcal K} }
\newcommand{\cL}{{\mathcal L} }
\newcommand{\cM}{{\mathcal M} }
\newcommand{\cN}{{\mathcal N} }
\newcommand{\cO}{{\mathcal O} }
\newcommand{\cP}{{\mathcal P} }
\newcommand{\cQ}{{\mathcal Q} }
\newcommand{\cS}{{\mathcal S} }
\newcommand{\cT}{{\mathcal T} }
\newcommand{\cW}{{\mathcal W} }
\newcommand{\cX}{{\mathcal X} }
\newcommand{\cY}{{\mathcal Y} }
\newcommand{\cZ}{{\mathcal Z} }
\newcommand{\bD}{{\bf D} }
\newcommand{\imp}{{\Rightarrow}}
\newcommand{\wt}{\widetilde}
\newcommand{\wh}{\widehat}
\newcommand{\Hom}{{\rm Hom}}
\def\ol#1{{\overline{#1}}}
\def\psh{{plurisubharmonic}}
\title{{\bf GLOBAL GENERATION OF THE DIRECT IMAGES OF RELATIVE PLURICANONICAL SYSTEMS}}
\date{December 4, 2010}
\author{Hajime TSUJI\footnote{Partially supported by Grant-in-Aid for Scientific Reserch (S) 17104001}}
\maketitle
\begin{abstract}
\noindent In this article, using the plurisubharmonic variation property of 
canonical measures (cf. \cite{canonical}), we prove 
that for an algebraic fiber space $f : X \to Y$, $f_{*}\mathcal{O}_{X}(mK_{X/Y})$ 
is globally generated on the complement of the discriminant locus of $f$ 
for every sufficiently large and divisible $m$.   As a byproduct, we prove Iitaka's conjecture on the subadditivity of Kodaira dimensions.      
\end{abstract}
\tableofcontents
\section{Introduction}

Let  $f : X \longrightarrow Y$ be a surjective projective morphism between  smooth projective varieties with connected fibers.   In this paper we shall call such a fiber space an {\bf algebraic fiber space} for simplicity.   We set $K_{X/Y} := K_{X}\otimes f^{*}K_{Y}^{-1}$
and call it the relative canonical line bundle of $f : X \to Y$.

Let $f : X \to Y$ be an algebraic fiber space.  
It is well known that the direct image $f_{*}\mathcal{O}_{X}(mK_{X/Y})$ is locally free outside of the discriminant locus (cf. \cite{s2,canAZD,can}) and is semipositive for every $m\geqq 1$ in certain algebraic senses (cf. Theorems \ref{kawamata} and \ref{viehweg} below). But Theorems \ref{kawamata} and \ref{viehweg}
 do not imply the existence of nontrivial global sections of 
   $f_{*}\mathcal{O}_{X}(mK_{X/Y})$. 
 
The purpose of this article is to prove that $f_{*}\mathcal{O}_{X}(mK_{X/Y})$ is globally generated  on the complement of the discriminant 
locus of $f$ for every sufficiently large and divisible $m$.

The main difficulty to prove the global generation is the fact that the direct 
image $f_{*}\mathcal{O}_{X}(mK_{X/Y})$ is only semipositive and not strictly positive (= ample) in general.   The idea of the proof is 
to distinguish the null direction of the positivity of $f_{*}\mathcal{O}_{X}(mK_{X/Y})$ as a Monge-Amp\`{e}re foliation and to realize 
the direct image $f_{*}\mathcal{O}_{X}(mK_{X/Y})$ (or its certain symmetric power)  as the pull back of an ample vector bundle on a certain moduli 
space via the moduli map.

\subsection{Kawamata's semipositivity theorem}

In order to clarify what is new in this article, I would like to review
 briefly  the former results and methods on 
the semipositivity of the direct images of pluricanonical systems in Sections 1.1 and 1.2.  

The first result on the semipositivity of the relative pluricanonical system 
is the following theorem due to Y. Kawamata in 1982.  

\begin{theorem}(\cite{ka1})\label{kawamata}
Let $f : X \to Y$ be an algebraic fiber space.   
Suppose that $\dim Y = 1$.   Then  for every positive integer $m$,  
$f_{*}\mathcal{O}_{X}(mK_{X/Y})$ is a semipositive vector bundle on $Y$, in the 
sense that every quotient $\mathcal{Q}$ of $f_{*}\mathcal{O}_{X}(mK_{X/Y})$, 
 $\deg \mathcal{Q} \geqq 0$ holds. 
\fbox{}
\end{theorem}
The proof of Theorem \ref{kawamata} depends on the variation of Hodge 
structure due to P.A. Griffiths  and W. Schmidt (cf. \cite{griff,sch}).
We note that before Theorem \ref{kawamata}, T. Fujita proved the case of 
$m=1$ in \cite{fu} by using the curvature computation of the  Hodge metrics of P.A. Griffiths (\cite{griff}). In this special case, Fujita gave a singular hermitian metric on the vector bundle $f_{*}\mathcal{O}_{X}(K_{X/Y})$ with semipositive 
curvature in the sense of Griffiths.   
In contrast to Fujita's result, for $m\geqq 2$, Theorem \ref{kawamata} does not give a (singular) hermitian metric on 
$f_{*}\mathcal{O}_{X}(mK_{X/Y})$ with semipositive curvature, because the proof 
relies on the semipositivity of the curvature of the Finslar metric 
on $f_{*}\mathcal{O}_{X}(mK_{X/Y})$ defined by  
\begin{equation}
\parallel\sigma\parallel := \left(\int_{X/Y}|\sigma|^{\frac{2}{m}}\right)^{\frac{m}{2}}   
\end{equation}
which is a singular hermitian metric on the tautological line bundle 
on $\mathbb{P}((f_{*}\mathcal{O}_{X}(mK_{X/Y}))^{*})$.
 
\subsection{Viehweg's semipositivity theorem}

In 1995 E. Viehweg extended Theorem \ref{kawamata} (\cite[Section 6]{v}) 
in the case of $f$-semiample relative canonical bundles   
 and constructed quasi-projective moduli spaces of polarized projective
 manifolds with semiample canonical bundles (\cite{v}).  
Since we use Viehweg's idea in this article, we  state his result 
precisely. First we  recall several definitions. 
\begin{definition}\label{globalgeneration}
Let $Y$ be a quasi-projective scheme, let $Y_{0}$ be an open dense suchscheme 
and let $\mathcal{G}$ be a coherent sheaf on $Y$, 
We say that $\mathcal{G}$ is {\bf globally generated} over $Y_{0}$,
if the natural map $H^{0}(Y,\mathcal{G})\otimes \mathcal{O}_{Y}
\to \mathcal{G}$ is surjective over $Y_{0}$. \fbox{}
\end{definition}
For a coherent sheaf $\mathcal{F}$ and a positive integer $a$, $S^{a}(\mathcal{F})$ denotes 
the $a$-th symmetric power of $\mathcal{F}$.  
To measure the positivity of coherent sheaves, we shall introduce the 
following notion.
\begin{definition}\label{weakpositivity}
Let $Y$ be a quasi-projective reduced scheme, $Y_{0}\subseteq Y$ an open dense subscheme and let $\mathcal{G}$ be locally free sheaf on $Y$, of finite constant  rank.  Then $\mathcal{G}$ is {\bf weakly positive} over $Y_{0}$, if 
for an ample invertible sheaf $\mathcal{H}$ on $Y$ and for a given number 
$\alpha > 0$ there exists some $\beta > 0$ such that 
$S^{\alpha\cdot\beta}(\mathcal{G})\otimes \mathcal{H}^{\beta}$ 
is globally generated over $Y_{0}$.  \fbox{}
\end{definition}  
The notion of weak positivity is a natural generalization of the notion of nefness of line bundles.  Roughly speaking, the weak semipositivity of $\mathcal{G}$ over $Y_{0}$ means that 
$\mathcal{G}\otimes \mathcal{H}^{\varepsilon}$ is $\mathbb{Q}$-globally generated over $Y_{0}$ for every $\varepsilon > 0$. 
\begin{definition}\label{succeq}
Let $\mathcal{F}$ be a locally free sheaf and let $\mathcal{A}$ be an invertible sheaf, both on a quasi-projective reduced scheme $Y$.  
We denote 
\begin{equation}
\mathcal{F} \succeq \frac{b}{a}\,\,\mathcal{A},
\end{equation}
if $S^{a}(\mathcal{F})\otimes \mathcal{A}^{-b}$ is weakly positive over $Y$,
where $a,b$ are positive integers.  
\fbox{}
\end{definition}
For a normal variety $X$, we define the canonical sheaf $\omega_{X}$ of 
$X$ by 
\begin{equation}
\omega_{X} := i_{*}\mathcal{O}_{X_{reg}}(K_{X_{reg}}), 
\end{equation}
where $X_{reg}$ denotes the regular part of $X$ and $i : X_{reg}\to X$ 
denotes the natural injection. 
The following notion introduced by Viehweg is closely related to the 
notion of logcanonical thresholds. 
\begin{definition}
Let $(X,\Gamma)$ be a pair of normal variety $X$ and an effective Cartier divisor $\Gamma$.  Let $\pi : X^{\prime}\to X$ be a log resolution of $(X,\Gamma)$ and let $\Gamma^{\prime}:= \pi^{*}\Gamma$. 
For a positive integer $N$ we define
\begin{equation}
\omega_{X}\left\{\frac{-\Gamma}{N}\right\}
= \pi_{*}\left(\omega_{X^{\prime}}\left(-\left\lfloor\frac{\Gamma^{\prime}}{N}\right\rfloor\right)\right)
\end{equation}
and 
\begin{equation}
\mathcal{C}_{X}(\Gamma,N)= \mbox{\em Coker}\left\{\omega_{X}\left\{\frac{-\Gamma}{N}\right\}\to \omega_{X}\right\}.
\end{equation}
If $X$ has at most rational singularities, one defines :
\begin{equation}
e(\Gamma)= \min\{N>0\,|\,\mathcal{C}_{X}(\Gamma,N) = 0\}.
\end{equation}
If $\mathcal{L}$ is an invertible sheaf,  $X$ is proper with at most rational  singularities and  $H^{0}(X,\mathcal{L})\neq 0$, then one defines
\begin{equation}
e(\mathcal{L}) = \sup\left\{ e(\Gamma)|\Gamma :\,\mbox{effective Cartier divisor 
with $\mathcal{O}_{X}(\Gamma)\simeq \mathcal{L}$}\right\}.
\end{equation}
\fbox{}
\end{definition}
Now we state the result of E. Viehweg. 
\begin{theorem}(\cite[p.191,Theorem 6.22]{v})\label{viehweg}
Let $f : X\to Y$ be a flat surjective projective Gorenstein morphism of reduced connected quai-projective schemes.  Assume that the sheaf $\omega_{X/Y}$
is $f$-semi-ample and that the fibers $X_{y} = f^{-1}(y)$ are reduced normal varieties with at most rational singularities.  Then one has :
\begin{enumerate}
\item[\em (1)]{\bf Functoriality}: For $m > 0$ the sheaf $f_{*}\omega_{X/Y}^{m}$ is locally free of 
rank $r(m)$ and it commutes with arbitrary base change.
\item[\em (2)]{\bf Weak semipositivity}: For $m > 0$ the sheaf $f_{*}\omega_{X/Y}^{m}$ is weakly positive over $Y$.
\item[\em (3)] {\bf Weak semistability}: Let $m > 1, e > 0$ and $\nu > 0$ be chosen so that 
$f_{*}\omega_{X/Y}^{m}\neq 0$ and 
\begin{equation}
e \geqq \sup\left\{\frac{k}{m-1}, e(\omega_{X_{y}}^{k})\,\,;\,\,\mbox{for}\,\,y\in Y\right\}
\end{equation}
hold. 
Then 
\begin{equation}
f_{*}\omega_{X/Y}^{m} \succeq \frac{1}{e\cdot r(k)}\,\det (f_{*}\omega_{X/Y}^{k})
\end{equation}
holds. \fbox{}
\end{enumerate}\vspace{3mm}
\end{theorem}

\noindent Although Theorem \ref{viehweg} assumes the $f$-semiampleness of $\omega_{X/Y}$, the advantages of this generalization are : 
\begin{itemize}
\item The base space is of arbitrary dimension.  
\item The semipositivity is more explicit than the one in Theorem \ref{kawamata}. 
\item The comparison of the positivity of $f_{*}\omega_{X/Y}^{m}$ 
and $\det (f_{*}\omega_{X/Y}^{m})$ is given.  
\end{itemize}

Later Theorems \ref{kawamata} and \ref{viehweg} have been extensively used in many 
other contexts (for example see \cite{subad,v}).
 
\subsection{Analytic Zariski decompositions}\label{AZD}

To state the main result, we  introduce
the notion of analytic Zariski decompositions.   
This notion will be used throughout this article. 

\begin{definition}\label{defAZD}
Let $M$ be a compact complex manifold and let $L$ be a holomorphic line bundle
on $M$.  A singular hermitian metric $h$ on $L$ is said to be 
an {\bf analytic Zariski decomposition}{\em (}{\bf AZD} in short{\em )}, if the followings hold.
\begin{enumerate}
\item[\em (1)] $\Theta_{h}$ is a closed positive current. 
\item[\em (2)] For every $m\geq 0$, the natural inclusion:
\begin{equation}
H^{0}(M,\mathcal{O}_{M}(mL)\otimes\mathcal{I}(h^{m}))\rightarrow
H^{0}(M,\mathcal{O}_{M}(mL))
\end{equation}
is an isomorphim. \fbox{}
\end{enumerate}
\end{definition}
\begin{remark} If an AZD exists on a line bundle $L$ on a smooth projective
variety $M$, $L$ is pseudoeffective by the condition 1 above. \fbox{}
\end{remark}

It is known that for every pseudoeffective line bundle on a compact complex manifold, there exists an AZD on $F$ (cf. \cite{tu,tu2,d-p-s}). 
The advantage of the AZD is that we can handle pseudoeffective line bundle 
$L$ on a compact complex manifold $X$  
as a singular hermitian  line bundle with semipositive curvature current
as long as we consider the ring $R(X,L) := \oplus_{m\geqq 0} H^{0}(X,\mathcal{O}_{X}(mL))$. 

We also note that there exists a smilar but different notion : {\em singular hemitian metrics with minimal singularities} introduced by Demailly, Peternell and Schneider (cf. \cite{d-p-s}).  A singular hemitian metric with minimal singularities is 
always an AZD, but in general an AZD need not be a singular hemitian metric with minimal singularities at least in the log canonical case (\cite{LC})  
\footnote{Actually this difference is closely related to the abundance conjecture.}.
In this article, we use the notion of AZD's, since 
the canonical measure (cf. Theorem \ref{canmeasure}) plays the crucial role in this article and the inverse of 
the canonical measure need not be a singular hermitian metric with minimal singularities.     
    
\subsection{Statement of the main results}
We note that Theorems \ref{kawamata} and \ref{viehweg} do not imply the existence of nontrivial global sections of 
$f_{*}\mathcal{O}_{X}(mK_{X/Y})$ for some $m > 0$.  
In this article we shall prove the global generation of $f_{*}\mathcal{O}_{X}(m!K_{X/Y})$ for every sufficiently large $m$ on the complement of the discriminant locus of $f$.    The following is the main result in this article.  
   
\begin{theorem}\label{main}
Let $f : X \longrightarrow Y$ be an algebraic fiber space and let
$Y^{\circ}$ be the complement of the discriminant locus of $f$ in $Y$.
Then we have the followings : 
\begin{enumerate}
\item[\em (1)]{\bf Global generation}:  
There exist positive integers $b$ and $m_{0}$ such that
for every integer $m$ satisfying $b\,|\,m$ and $m \geqq m_{0}$,   
$f_{*}\mathcal{O}_{X}(mK_{X/Y})$ is globally generated over  
$Y^{\circ}$.
\item[\em (2)]{\bf Weak semistability 1}:  Let  $m$ be a positive integer such that $f_{*}\mathcal{O}_{X}(mK_{X/Y})\neq 0$.  Let $r$ denote $\mbox{\em rank}\,f_{*}\mathcal{O}_{X}(mK_{X/Y})$
and let $X^{r}:= X\times_{Y}X\times_{Y}\cdots\times_{Y}X$ be the $r$-times fiber product over $Y$.  Let $f^{r}: X^{r}\to Y$ be the natural morphism. 

Let $\Gamma \in |mK_{X^{r}/Y}-f^{r*}\det f_{*}\mathcal{O}_{X}(mK_{X/Y})|$
be the effective divisor corresponding to the canonical inclusion :
\begin{equation}\label{wp1}
f^{r*}(\det f_{*}\mathcal{O}_{X}(mK_{X/Y}))\hookrightarrow f^{r*}f^{r}_{*}\mathcal{O}_{X^{r}}(mK_{X^{r}/Y}) \hookrightarrow \mathcal{O}_{X^{r}}(mK_{X^{r}/Y}).
\end{equation}
Then $\Gamma$ does not contain any fiber $X^{r}_{y} (y\in Y^{\circ})$ such that if we we define the number $\delta_{0}$ by 
\begin{equation}
\delta_{0} := \sup \{\delta \,|\,(X^{r}_{y},\delta\cdot\Gamma_{y}) \,\,\mbox{is KLT for all $y\in Y^{\circ}$}\}, 
\end{equation}
then for every $\varepsilon < \delta_{0}$ and a sufficiently large positive integer $d$,  
\begin{equation}\label{wp2}
f_{*}\mathcal{O}_{X}(d!K_{X/Y})\succeq 
\frac{d!\varepsilon}{(1 + m\varepsilon)r}\,\det f_{*}\mathcal{O}_{X}(mK_{X/Y})
\end{equation}
holds over $Y^{\circ}$. 
\item[\em (3)]{\bf Weak semistability 2}: There exists a singular hermitian metric $H_{m,\varepsilon}$ on 
$(1+m\varepsilon)K_{X^{r}/Y} - \varepsilon\cdot f^{r*}\det f_{*}\mathcal{O}_{X}(mK_{X/Y})^{**}$
such that 
\begin{enumerate}
\item[\em (a)] $\sqrt{-1}\,\Theta_{H_{m,\varepsilon}}\geqq 0$ holds on $X^{r}$ in the sense 
of current.
\item[\em (b)] For every $y\in Y^{\circ}$, $H_{m,\varepsilon}|_{X^{r}_{y}}$ is well defined 
and is an AZD (cf. Definition \ref{AZD}) of 
\begin{equation}
(1+m\varepsilon)K_{X^{r}/Y} - \varepsilon\cdot f^{r*}\det f_{*}\mathcal{O}_{X}(mK_{X/Y})^{**}|X_{y}. 
\end{equation}
\end{enumerate} 
\fbox{} 
\end{enumerate}
\end{theorem}
\begin{remark}
The 3rd assertion implies the 2nd assertion. \fbox{} 
\end{remark}
The major difference between Theorems \ref{main} and  \ref{kawamata} is 
that  in Theorem \ref{main} $f_{*}\mathcal{O}_{X}(mK_{X/Y})$ is globally generated over the complement of the discriminant locus of $f$, while Theorem \ref{kawamata} implies the semipositivity of $f_{*}\mathcal{O}_{X}(mK_{X/Y})$.
In this sense Theorem \ref{main} is much stronger than Theorem \ref{kawamata}. 
The major difference between Theorems \ref{main} and \ref{viehweg} is 
(besides the global generation assertion) that  in Theorem \ref{main}, 
we do not assume the $f$-semiampleness of $K_{X/Y}$ in Theorem \ref{main}. 
\vspace{3mm} \\  
\noindent We also have the following log version of Theorem \ref{main}. 

\begin{theorem}\label{logmain}
Let $f : X \to Y$ be an algebraic fiber space and let
$D$ be an effective $\mathbb{Q}$ divisor on $X$ 
such that $(X,D)$ is KLT.  Let $Y^{\circ}$ denote the complement of the 
discriminant locus of $f$. 
We set  
\begin{equation}
Y_{0}:= \{y\in Y|y\in Y^{\circ}, (X_{y},D_{y})\,\,\mbox{is a KLT pair}\} 
\end{equation}
\begin{enumerate}
\item[\em (1)]{\bf Global generation}:  
There exist positive integers $b$ and $m_{0}$ such that 
for every for every integer $m$ satisfying $b\,|\,m$ and  $m\geqq m_{0}$, $m(K_{X/Y}+D)$ is Cartier and 
$f_{*}\mathcal{O}_{X}(m(K_{X/Y}+D))$ is globally generated 
over $Y_{0}$.    
\item[\em (2)]{\bf Weak semistability 1}: 
Let $m$ be a positive integer such that $m(K_{X/Y}+D)$ is integral and
$f_{*}\mathcal{O}_{X}(m(K_{X/Y}+D))\neq 0$.
  Let $r$ denote $\mbox{\em rank}\,f_{*}\mathcal{O}_{X}(\lfloor m(K_{X/Y}+D)\rfloor)$. 
Let $X^{r}:= X\times_{Y}X\times_{Y}\cdots\times_{Y}X$ be the $r$-times fiber product over $Y$ and let $f^{r}: X^{r}\to Y$ be the natural morphism. 
And let $D^{r}$ denote the divior on $X^{r}$ defined by 
$D^{r} = \sum_{i=1}^{r} \pi_{i}^{*}D$,  
where $\pi_{i} : X^{r} \longrightarrow X$ denotes the projection: 
$X^{r}\ni (x_{1},\cdots ,x_{n}) \mapsto x_{i} \in X$.  

There exists a canonically defined effective divisor $\Gamma$ (depending on $m$) on $X^{r}$ which does not conatin any fiber $X^{r}_{y} (y\in Y^{\circ})$ such that if we we define the number $\delta_{0}$ by 
\begin{equation}
\delta_{0} := \sup \{\delta\, |\, (X^{r}_{y},D^{r}_{y}+\delta\Gamma_{y})\,\, \mbox{is KLT for all $y\in Y^{\circ}$}\}, 
\end{equation}
then for every $\varepsilon < \delta_{0}$ and every sufficiently large positive integer $d$,  
\begin{equation}
f_{*}\mathcal{O}_{X}( d!(K_{X/Y}+D))\succeq \frac{d!\varepsilon}{(1+m\varepsilon)r}\,\det f_{*}\mathcal{O}_{X}(\lfloor m(K_{X/Y}+D)\rfloor )
\end{equation}
holds over $Y_{0}$.
\item[\em (3)]{\bf Weak semistability 2}:  There exists a singular hermitian metric $H_{m,\varepsilon}$ on 
\begin{equation}
(1+m\varepsilon)(K_{X^{r}/Y}+D^{r}) - \varepsilon\cdot f^{*}\det f_{*}\mathcal{O}_{X}(\lfloor m(K_{X/Y}+D)\rfloor)^{**}
\end{equation}
such that 
\begin{enumerate}
\item $\sqrt{-1}\,\Theta_{H_{m,\varepsilon}}\geqq 0$ holds on $X$ in the sense 
of current.
\item  
 For every $y\in Y_{0}$, $H_{m,\varepsilon}|X^{r}_{y}$ is well defined 
and is an AZD of
\[ 
(1+m\varepsilon)(K_{X^{r}/Y}+D^{r}) - \varepsilon\cdot f^{r*}\det f_{*}\mathcal{O}_{X}(\lfloor m(K_{X/Y}+D)\rfloor)^{**}|X_{y}
\]
\fbox{} 
\end{enumerate} 

\end{enumerate}
\end{theorem}

\noindent The main ingredient  of the proof of Theorems \ref{main} and \ref{logmain} is the (logarithmic) plurisubharmonic variation property of canonical measures 
(Theorem \ref{relative} in \cite{canonical}). 
The new feature of the proof is the use of the Monge-Amp\`{e}re foliations 
arising from the canonical measures and the weak semistability of the 
direct images of relative pluricanonical systems.   
One may consider these new tools as substitutes of the local Torelli theorem 
for  minimal models with semiample canonical divisors in \cite{ka3}. \vspace{3mm}\\ 

The scheme of the proof is as follows.   For an algebraic fiber space 
$f : X \to Y$ with $\mbox{Kod}(X/Y)\geqq 0$ (cf. (\ref{rel})), we take the relative 
canonical measure $d\mu_{can,X/Y}$ (see Section \ref{relativemeasure}).  Then the null distribution of the curvature 
$\Theta_{d\mu_{can,X/Y}^{-1}}$ of the singular hermitian metric 
$d\mu_{can,X/Y}^{-1}$ on $K_{X/Y}$ defines a singular Monge-Amp\`{e}re foliation 
on $X$.   Here the important fact  is that the leaf of the foliation is 
complex analytic (\cite{b-k}) (although it is not clear that the foliation itself is complex analytic apriori). By using the weak semistability of $f_{*}\mathcal{O}_{X}(m!K_{X/Y})$, 
we  prove that this singular foliation actually descends to a 
singular foliation $\mathcal{G}$ on the base space $Y$.  
Let us define the  (singular) hermitian metric $h_{m}$ on 
 $f_{*}\mathcal{O}_{X}(m!K_{X/Y})$ defined by 
 \begin{equation}
h_{m}(\sigma,\sigma^{\prime}) := \int_{X/Y} \sigma\cdot\overline{\sigma^{\prime}}\cdot d\mu_{can,X/Y}^{-(m!-1)}. 
\end{equation}
Then we see that $(f_{*}\mathcal{O}_{X}(m!K_{X/Y}),h_{m})$ is flat 
along the leaves of $\mathcal{G}$ on $Y$. 
Taking  $m$ sufficiently large, we see that 
the metrized relative canonical model (cf. Definition \ref{metrizedcanonical} below) of $f : X \to Y$ is locally trivial along 
the leaves. Then we see that the leaves of $\mathcal{G}$ consists of the fiber of the moduli map to the moduli space of relative canonical models marked 
with the metrized Hodge line bundles.  Then the global generation property 
of $f_{*}\mathcal{O}_{X}(mK_{X/Y})$ follows from the Nakai-Moishezon type argument. 

\subsection{Iitaka's conjecture}

In this subsection, we apply Theorem \ref{main} to Iitaka's conjecture.
The following conjecture by S. Iitaka (\cite{i}) is  well known. 
 
\begin{conjecture}(Iitaka's conjecture)\label{iitaka}
Let $f : X \longrightarrow Y$ be an algebraic fiber space.  
Then 
\begin{equation}
\mbox{\em Kod}(X) \geqq \mbox{\em Kod}(Y) +\mbox{\em Kod}(X/Y)
\end{equation}
holds, where $\mbox{\em Kod}(X),\mbox{\em Kod}(Y)$
denote the Kodaira dimension (cf. (\ref{kod})) of $X,Y$ repsectively and 
$\mbox{\em Kod}(X/Y)$ denotes the relative Kodaira dimension as  (\ref{rel}). \fbox{}
\end{conjecture} 
The typical examples of algebraic fiber spaces are Iitaka fibrations, Albanese 
maps, the universal families over fine moduli spaces. 
Especially the Iitaka fibration $f : X \to Y$ has the property that 
$\mbox{Kod}(X/Y) = 0$.  Hence Conjecture \ref{iitaka} reduces the birational 
classification of $X$ to the study of families of varieties with Kodaira dimension $0$ and the study of the base sace $Y$ with $\mbox{Kod}(Y)\leqq \mbox{Kod}(X)$.   
Conjecture \ref{iitaka} is considered to be one of the key 
for the birational classification of projective varieties . 
For detailed explanation and references, see the survey article \cite{mo} for example.

In \cite{ka1} Kawamata solved Conjecture \ref{iitaka} 
in the case of $\dim Y = 1$ by using Theorem \ref{kawamata}.  
And if $\mbox{Kod}(Y) = \dim Y$,i.e., $Y$ is of general type, 
then Conjecture \ref{iitaka} can be easily deduced from 
Theorem \ref{kawamata}.  And in the case that  a general fiber of $f : X \to Y$ is of general type, Conjecture \ref{iitaka} has been solved (cf. \cite{v1,ko}). 
 And in \cite{ka3}, Kawamata reduced Conjecture 
\ref{iitaka} to the completion of the minimal model program (MMP).
Hence by the completion of MMP in dimension 3 (see \cite{k-m} for example),  Conjecture \ref{iitaka} 
has been solved in the case of $\dim X = 3$.  

As an immediate consequence of Theorem \ref{main}, we give an   
affirmative answer to Iitaka's conjecture. 
\begin{theorem}\label{iitakas}
Conjecture \ref{iitaka} holds.  \fbox{}
\end{theorem}
{\em Proof of Theorem \ref{iitakas}}.
Let $f : X \to Y$ be an algebraic fiber space.  
If $\mbox{Kod}(Y)$ or $\mbox{Kod}(X/Y)$ is $-\infty$, Conjecture \ref{iitaka} 
certainly holds.  Hence we  assume that $\mbox{Kod}(Y)$ and $\mbox{Kod}(X/Y)$ are nonnnegative. 
We note that there exists a natural morphism :  
\begin{equation}
H^{0}(Y,f_{*}\mathcal{O}_{X}(mK_{X/Y}))
\otimes 
H^{0}(Y,\mathcal{O}_{Y}(mK_{Y})) \to 
H^{0}(X,\mathcal{O}_{X}(mK_{X})). 
\end{equation}
Then by Theorem \ref{main}, we have that  
\begin{equation}
\limsup_{m\to\infty}\frac{\log \dim H^{0}(Y,f_{*}\mathcal{O}_{X}(mK_{X/Y}))}
{\log m} \geqq \mbox{Kod}(X/Y)
\end{equation}
holds. 
Hence we see that 
\begin{equation}
\mbox{Kod}(X) \geqq \mbox{Kod}(Y) + \mbox{Kod}(X/Y)
\end{equation}
holds. \fbox{}
\begin{remark}
The optimal form of Iitaka's conjecture is:
\begin{equation}
\mbox{\em Kod}(X) \geqq \mbox{\em Kod}(Y) +\max\{\mbox{\em Kod}(X/Y),\mbox{\em Var}(f)\}.
\end{equation}
At this moment, I do not know the proof.  
\fbox{}
\end{remark}

The organization of this article is as follows.  In Section 2, we review the 
canonical measures intorduced in \cite{s-t,canonical}.  Especially the 
logarithmic subhrmonicity of the canonical meaures (cf. \cite{canonical,LC}) is  explained.   Using the logarithmic subharmonicity and Viehweg's idea, we prove
the weak semistability of the direct images of relative pluri log canonical systems for a family of KLT pairs.  
In Section 3, we construct the moduli space of the metrized canonical models of KLT pairs.  The construction is rather standard, but technical. 
In Section 4, we analyse the Monge-Amp\`{e}re foliation assuming the regularity
results of canonical measures which is proven Section 6 below. 
In Section 5, we complete the proof of the main results assuming the regularity  rusults in Section 6. 
In Section 6, we prove the regularity of canonical measures by the 
dynamical construction of canonical measures and H\"{o}rmander's $L^{2}$-estimate of $\bar{\partial}$-operators.  
In Section 7, we provide several technical results which are used in 
Section 4.\\

\noindent The order of contents may be a little bit irregular.  But 
I hope that to put off the technical stuffs later makes the scheme of the proof clear. \vspace{3mm} \\

\noindent{\bf Notations}
\begin{itemize}
\item  For a real number $a$, $\lceil a\rceil$ denotes the minimal integer greater than or equal to $a$ and $\lfloor a\rfloor$ denotes the maximal integer smaller than 
or equal to $a$.  
\item Let $X$ be a projective variety and let $D$ be a Weil divsor on $X$.
Let $D = \sum d_{i}D_{i}$ be the irreducible decomposition.  
We set 
\begin{equation}
\lceil D\rceil := \sum \lceil d_{i}\rceil D_{i}
,\lfloor D \rfloor := \sum \lfloor d_{i}\rfloor D_{i}.
\end{equation} 
\item Let $f : X \to Y$ be an algebraic fiber space and let $D$ be 
a $\mathbb{Q}$-divisor on $X$.  Let 
\begin{equation}\label{h-v}
D = D^{h} + D^{v}
\end{equation}
be  the decomposition such that an irreducible component of 
$\mbox{Supp}\, D$ is contained in $\mbox{Supp}\, D^{h}$ if and only if 
it is mapped onto $Y$.  $D^{h}$ is the horizontal part of $D$ 
and $D^{v}$ is the vertical part of $D$.  
\item Let $(X,D)$ be a pair of a normal variety and a $\mathbb{Q}$-divisor 
on $X$.  Suppose that $K_{X} + D$ is $\mathbb{Q}$-Cartier. 
Let $f : Y \to X$ be a log resolution.  Then we have the formula :
\[
K_{Y} = f^{*}(K_{X}+D) + \sum a_{i}E_{i}, 
\]
where $E_{i}$ is a prime divisor and $a_{i}\in \mathbb{Q}$. 
The pair $(X,D)$ is said to be {\bf subKLT}(resp. {\bf subLC}), if $a_{i} > -1$
(resp. $a_{i} \geqq -1$) holds for every $i$. 
$(X,D)$ is said to be {\bf KLT} (resp. {\bf LC}), if $(X,D)$ is subKLT(resp. subLC) and $D$ is effective.
\item Let $X$ be a projective variety and let $\mathcal{L}$ be an invertible
sheaf on $X$. $\mathcal{L}$ is said to be semiample, if there exists a positive
integer $m$ such that $|\mathcal{L}^{\otimes m}|$ is  base point free.  
\item $f : X \to Y$ be a morphism between projective varieties. 
Let $\mathcal{L}$ be an invertible sheaf on $X$.  $\mathcal{L}$ is said to 
be $f$-semiample, if for every $y \in Y$, $\mathcal{L}|f^{-1}(y)$ is 
semiample.  
\item Let $L$ be a $\mathbb{Q}$-line bundle on a compact complex manifold $X$, i.e., $L$ is a formal fractional power of a genuine line bundle on $X$.  
A  singular hermitian metric $h$ on $L$ is given by
\[
h = e^{-\varphi}\cdot h_{0},
\]
where $h_{0}$ is a $C^{\infty}$ hermitian metric on $L$ and 
$\varphi\in L^{1}_{loc}(X)$ is an arbitrary function on $X$.
We call $\varphi$ the  weight function of $h$ with respect to $h_{0}$.  We note that $h$ makes sense,
 since a hermitian metric is a real object. 

The curvature current $\Theta_{h}$ of the singular hermitian $\mathbb{Q}$-line
bundle $(L,h)$ is defined by
\[
\Theta_{h} := \Theta_{h_{0}} + \partial\bar{\partial}\varphi,
\]
where $\partial\bar{\partial}\varphi$ is taken in the sense of current.
We define the multiplier ideal sheaf  ${\cal I}(h)$ of $(L,h)$  
by 
\[
\mathcal{I}(h)(U) := \{ f \in \mathcal{O}_{X}(U); \,\, |f|^{2}\,e^{-\varphi}
\in L^{1}_{loc}(U)\},
\]
where $U$ runs open subsets of $X$.
\item A singular hermitian line bundle $(L,h)$ is said to be {\bf pseudoeffective}, if $\sqrt{-1}\Theta_{h}$ is a closed semipositive current.  
\item For a closed positive $(1,1)$ current $T$, $T_{abc}$ denotes 
the abosolutely continuous part of $T$. 
\item  For  a Cartier divisor $D$, we denote the 
corresponding line bundle by the same notation. 
Let $D$ be an effective $\mathbb{Q}$-divisor on a smooth projective variety $X$.  Let $a$ be a positive integer such that $aD$ is Cartier. 
We identify $D$ with a formal $a$-th root of the line bundle $aD$.  
We say that $\sigma$ is a multivalued global holomorphic section of $D$ 
with divisor $D$, if  $\sigma$ is a formal $a$-th root of a global holomorphic section of $aD$ with divisor $aD$. And $1/|\sigma|^{2}$ denotes the singular hermitian metric on $D$ defined by 
\[
\frac{1}{|\sigma|^{2}} := \frac{h_{D}}{h_{D}(\sigma,\sigma)}, 
\]
where $h_{D}$ is an arbitrary $C^{\infty}$ hermitian metric on $D$.
\item For a singular hermitian line bundle $(F,h_{F})$ on a compact complex 
manifold $X$ of dimension $n$.   $K(X,K_{X}+F,h_{F})$ denotes (the diagonal part of) the Bergman kernel of $H^{0}(X,{\cal O}_{X}(K_{X} + F)\otimes {\cal I}(h_{F}))$ 
with respect to the $L^{2}$-inner product: 
\begin{equation}\label{inner}
(\sigma ,\sigma^{\prime}) := (\sqrt{-1})^{n^{2}}\int_{X}h_{F}\cdot\sigma\wedge \bar{\sigma}^{\prime}, 
\end{equation}  
i.e., 
\begin{equation}\label{BergmanK}
K(X,K_{X}+F,h_{F}) = \sum_{i=0}^{N}|\sigma_{i}|^{2}, 
\end{equation}
where $\{\sigma_{0},\cdots ,\sigma_{N}\}$ is a complete orthonormal basis 
of $H^{0}(X,{\cal O}_{X}(K_{X} + F)\otimes {\cal I}(h_{F}))$. 
It is clear that $K(X,K_{X}+F,h_{F})$ is independent of the choice of 
the complete orthonormal basis. 
\end{itemize}
     
\section{Canonical measures}
In this section we  review the definition 
and the basic properties of canonical measures
which plays the key role\footnote{Probably we may use the Narashimhan-Simha volume form (\cite{n-s}) instead of canonical measures to prove  Theorem \ref{main} and \ref{logmain}.  
For a smooth projective vairiety $X$ with nonnegative Kodaira dimension
and a positive integer $m$, the $m$-th Narashimhan-Simha volume form $K_{m}^{NS}$ is defined by  
\begin{equation}
K_{m}^{NS}(x) := \{|\sigma|^{\frac{2}{m}}(x)|  
\sigma \in H^{0}(X,\mathcal{O}_{X}(mK_{X})), \int_{X}|\sigma|^{\frac{2}{m}} =1\}. 
\end{equation}
The advantage of the  Narashimhan-Simha volume form is that 
its construction is much simpler than the one of the canonical measure.   But on the other hand, it seems to be hard to prove the 
regularity of the Narashimhan-Simha measure.}  
 in the proof of Theorems \ref{main} and \ref{logmain}.

The canonical measure is a natural generalization of K\"{a}hler-Einstein volume form to the case of projective varieties with nonnegative Kodaira dimension
(cf. \cite{s-t,canonical}).
The basic properties of the canonical measure are : 
\begin{enumerate}
\item[(1)] It is completely determined by the complex structure of the variety 
and is birationally invariant.  
\item[(2)] It is $C^{\infty}$ on a on a nonempty Zariski open subset 
of the variety and satisfies a Monge-Amp\`{e}re equation on a Zariski open 
subset on the base space of the Iitaka fibration (cf. Section \ref{IF}). 
\item[(3)] The logarithm of the measure is plurisubharmonic under projective 
deformations. 
\end{enumerate} 
For the detailed account, see \cite{s-t,canonical,LC}. 
The canonical measure is defined on 
an arbitrary KLT pair with nonnegative logarithmic Kodaira dimension 
(\cite{LC}).

\subsection{Iitaka fibration}\label{IF}

To construct the canonical measure, we need to  consider the Iitaka fibration.  The Iitaka fibration is the most naive way to extract the positivity of 
the canonical bundle on  a smooth projective variety with nonnegative 
Kodaira dimension. 

Let $X$ be a smooth projective variety. 
We define the Kodaira dimension of $X$ by 
\begin{equation}\label{kod}
\mbox{Kod}(X):= \limsup_{m\rightarrow\infty}\frac{\log \dim H^{0}(X,\mathcal{O}_{X}(mK_{X}))}{\log m}. 
\end{equation}
More generally for a KLT pair $(X,D)$, we define {\em the Kodaira dimension} 
$\mbox{Kod}(X,D)$ of $(X,D)$ by 
\begin{equation}
\mbox{Kod}(X,D):= \limsup_{m\rightarrow\infty}\frac{\log \dim H^{0}\left(X,\mathcal{O}_{X}(\lfloor m(K_{X}+D)\rfloor)\right)}{\log m}.
\end{equation}
Similarly for an algebraic fiber space $f : X \to Y$, we define 
{\em the relative Kodaira dimension} $\mbox{Kod}(X/Y)$ by 
\begin{equation}\label{rel}
\mbox{Kod}(X/Y) := \mbox{Kod}(F),
\end{equation}
where $F$ is a general fiber of $f$. 

Let $X$ be a smooth projective variety with $\mbox{Kod}(X)\geqq 0$. 
Then for a sufficiently large $m > 0$, the complete linear system 
$|m!K_{X}|$ gives a rational fibration (with connected fibers) :
\begin{equation}
f : X -\cdots\rightarrow  Y. 
\end{equation}
We call $f : X -\cdots\rightarrow  Y$ the {\bf Iitaka fibration} of $X$.
  
The Iitaka fibration is independent of the choice of the sufficiently large $m$
up to birational equivalence.  See \cite{i} for detail.   
In this sense the Iitaka fibration is unique. 
By taking a suitable modification, we may assume that 
$f$ is a morphism and $Y$ is smooth.   

The Iitaka fibration $f : X \to Y$ satisfies the following properties: 
\begin{enumerate}
\item[(1)] For a general fiber $F$, $\mbox{Kod}(F) = 0$ holds,    
\item[(2)] $\dim Y = \mbox{Kod}(X)$. 
\end{enumerate}

\subsection{Relative Iitaka fibrations}

The Iitaka fibration can be easily generalized to the relative setting.
This generalization will be used to analyze the variation of canonical measures
 on a projective faimily.  

Let $f : X \longrightarrow Y$ be an algebraic fiber space, i.e., 
$X,Y$ are smooth projective varieties  and $f$ is a proper surjective morphism 
with connected fibers.  

Let $m$ be a sufficiently large positive integer and we set  $F_{m}:= f_{*}\mathcal{O}_{X}(m!K_{X/Y})^{**}$. 
For $x\in X$, we set  
\begin{equation}
ev_{x} : F_{m,f(x)} \longrightarrow mK_{X/Y} 
\end{equation}
be the evaluation map. 
We define the relative canonical map : 
\begin{equation}
g : X -\cdots \rightarrow  \mathbb{P}(F_{m}^{*})
\end{equation} 
by 
\begin{equation}
g(x) := \{ [v^{*}]|\, v^{*}\in  F_{m,f(x)}^{*},\, v^{*}\!|\mbox{Ker}\,\,ev_{x} = 0\}.
\end{equation}
Let $Z$ be the image of $g$.  
Then we have the commutative diagram :  
\begin{equation}\label{triangle}
\begin{picture}(400,400)
\setsqparms[1`1`1`1;350`350]
\putVtriangle(0,00)%
[X` Z ` Y;g `f `h]
\end{picture}
\end{equation}  
For a sufficiently large $m$, we see that a general fiber $F$ of 
$g : X -\cdots\rightarrow Z$ is connected and $\mbox{Kod}(F) = 0$. 
We call $g : X -\cdots\rightarrow Z$@{\bf  a@relative Iitaka fibration}. 
By taking a suitable modification of $X$, we may assume that 
$g$ is a morphism. 

Let $f : X \longrightarrow Y$ be an algebraic fiber space 
and let $g : X \longrightarrow Z$ be a relative Iitaka fibration 
associated with $f_{*}\mathcal{O}_{X}(m!K_{X/Y})$. 
Taking a suitable modification we may and do assume the 
followings :
\begin{enumerate}
\item[(1)] $g$ is a morphism. 
\item[(2)] $Z$ is smooth. 
\item[(3)] $g_{*}\mathcal{O}_{X}(m!K_{X/Z})^{**}$ is a line bundle on $Z$
for every sufficiently large $m$. 
\end{enumerate}  
Let $h : Z \longrightarrow Y$ be the natural morphism.

This construction can be easily generalized to the case of a  
KLT pair $(X,D)$ with algebraic fiber space structure 
$f : X \to Y$. 

Also by the finite generation of the log canonical ring for a  KLT pair 
(\cite{b-c-h-m}), we may take $Z$ be be a family of logcanonical models
at least on a nonempty Zariski open subset of $Y$.  
In this case $Z$ has singularities. 
We call such a triangle (\ref{triangle}) or the family $h: Z \to Y$, {\bf the  relative canonical model}. 

\subsection{Hodge line bundles associated with Iitaka fibrations}
\label{hodgebundle}

Let $f : X \to Y$ be an  Iitaka fibration such that $X,Y$ are smooth 
and $f$ is a morphism. 
Then by \cite[p.169,Proposition 2.2]{f-m}, 
$f_{*}\mathcal{O}_{X}(m!K_{X/Y})^{**}$  is locally free on $Y$ for every  sufficiently large $m$, where $**$ denotes the double dual.
Since $f : X \to Y$ is an Iitaka fibration, a general fiber is of Kodaira 
dimension $0$ and  the direct image 
$f_{*}\mathcal{O}_{X}(m!K_{X/Y})$ is of rank $1$ for every sufficiently large
$m$.  
We define the $\mathbb{Q}$-line bundle $L_{X/Y}$ on $Y$ by  
\begin{equation}\label{HL}
L_{X/Y} := \frac{1}{m!}\,f_{*}\mathcal{O}_{X}(m!K_{X/Y})^{**}.  
\end{equation}
We note that $L_{X/Y}$ is independent of a sufficiently large $m$ (cf. \cite[Section 2]{f-m}). Let us fix such a $m$. 
Let $Y^{\circ}$ denote the complement of the discriminant locus of 
$f : X \to Y$. 
Then $L_{X/Y}$  carries the natural singular hermitian metric $h_{L_{X/Y}}$ defined by 
\begin{equation}\label{HM}
h_{L_{X/Y}}^{m!}(\sigma,\sigma)_{y} := \left(\int_{X_{y}}|\sigma|^{\frac{2}{m!}}\right)^{m!}, 
\end{equation}
where $y \in Y^{\circ}, X_{y}:= f^{-1}(y)$ and $\sigma \in m!L_{X/Y,y}$.  $h_{L_{X/Y}}$ is defined on 
$L_{X/Y}|Y^{\circ}$ apriori.  But by the theory of variation of Hodge structures (\cite{sch}), $h_{L_{X/Y}}$ extends to a singular hermitian metric 
on $L_{X/Y}$.   
It is known that $h_{L_{X/Y}}$ has semipositive curvature in the sense of current
(\cite{ka1}).
 
\subsection{Definition of canonical measures and the existence}
Now we define the canonical semipositive current
on a smooth projective variety of nonnegative Kodaira dimension. 
Let $f : X \to Y$ be the Iitaka fibration such that some positive multiple of the Hodge $\mathbb{Q}$-line bundle
$L_{X/Y}$ defined as in the last subsection is locally free. 

\begin{theorem}\label{canmeasure}(cf. \cite[Theorem 1.5]{canonical} and \cite[Theorem B.2]{s-t}) 
In the above notations, there exists a unique singular hermitian metric 
on $h_{K}$ on $K_{Y}+L_{X/Y}$ such that 
\begin{enumerate}
\item[\em (1)] $h_{K}$ is an AZD of $K_{Y} + L_{X/Y}$,
\item[\em (2)] $f^{*}h_{K}$ is an AZD of $K_{X}$,   
\item[\em (3)] $h_{K}$ is $C^{\infty}$ on a nonempty Zariski open subset $U$,
\item[\em (4)] $\omega_{Y} = \sqrt{-1}\,\Theta_{h_{K}}$ is a K\"{a}hler form 
on $U$,
\item[\em (5)] $-\mbox{\em Ric}_{\omega_{Y}} + \sqrt{-1}\,\Theta_{h_{L_{X/Y}}} = \omega_{Y}$
holds on $U$, where $h_{L_{X/Y}}$ denotes the Hodge metric defined as (\ref{HM}). \fbox{} 
\end{enumerate} \vspace{3mm}
\end{theorem}

\noindent The above equation:
\begin{equation}\label{LKE}
-\mbox{Ric}_{\omega_{Y}} + \sqrt{-1}\,\Theta_{h_{L_{X/Y}}} = \omega_{Y}
\end{equation}
is similar to the K\"{a}hler-Einstein equation :
\begin{equation}
-\mbox{Ric}_{\omega_{Y}}  = \omega_{Y}. 
\end{equation}
The correction term $\sqrt{-1}\,\Theta_{h_{L_{X/Y}}}$ reflects the 
isomorphism :
\begin{equation}\label{caniso}
R(X,K_{X})^{(a)} = R(Y,K_{Y}+L_{X/Y})^{(a)} 
\end{equation}
for some positive integer $a$, where for a graded ring $R : = \oplus_{i=0}^{\infty}R_{i}$
 and a positive integer $b$,  we set 
\begin{equation}
R^{(b)} : = \oplus_{i=0}^{\infty}R_{bi}.
\end{equation}

Now we shall define the canonical measure. 

\begin{definition}(\cite{s-t,canonical,LC})\label{L-K-E}
The current $\omega_{Y}$ on $Y$ constructed in Theorem \ref{canmeasure} is said to be {\bf the  canonical K\"{a}hler current}  of the Iitaka fibration $f : X \longrightarrow Y$.  Also $\omega_{X}:= f^{*}\omega_{Y}$ is said to be {\bf the canonical semipositive current} on $X$. 
We define the measure $d\mu_{can}$ on $X$ by 
\begin{equation}
d\mu_{can}:= \frac{1}{n!}f^{*}\left(\omega_{Y}^{n}\cdot h_{L_{X/Y}}^{-1}\right)
\end{equation}
and is said to be {\bf the canonical measure}, where $n$ denotes 
$\dim Y$. Here we note that $\omega_{Y}^{n}$ is a degenerate volume form on $Y$ and  
$f^{*}h_{X/Y}^{-1}$ is considered to be a relative (degenerate) volume form on 
$f : X \to Y$ (cf. (\ref{HM})), hence  $f^{*}\left(\omega_{Y}^{n}\cdot h_{L_{X/Y}}^{-1}\right)$ 
is considered to be a degenerate volume form on $X$. 
\fbox{} 
\end{definition}

We also have the log version of Theorem \ref{canmeasure} which plays a 
crucial role  not only  in the proof of Theorem \ref{logmain}, but also 
in the one of Theorem \ref{main}. 

Let $(X,D)$ be a KLT pair such that $X$ is smooth projective. 
We assume that $\mbox{Kod}(X,D) \geqq 0$, i.e., for every $m >> 1$, 
$|m!(K_{X}+D)| \neq \emptyset$.  
Let  
\begin{equation}
f : X \longrightarrow Y
\end{equation}
be a log Iitaka fibration of $(X,D)$.    After modifications, we may assume
the followings: 
\begin{enumerate}
\item[\em (1)] $X$,$Y$ are smooth and $f$ is a morphism with connected fibers. 
\item[\em (2)] $\mbox{Supp}\,D$ is a divisor with normal crossings. 
\item[\em (3)] There exists an effective divisor $\Sigma$ on $Y$ such that 
$f$ is smooth over $Y - \Sigma$, $\mbox{Supp}\, D^{h}$ is relatively normal crossings
over $Y - \Sigma$ and $f(D^{v})\subset \Sigma$, where $D^{h},D^{v}$ denote the horizontal and the vertical component of $D$ respectively.
\item[\em (4)] There exists a positive integer $m_{0}$ such that $f_{*}\mathcal{O}_{X}\left(m!(K_{X/Y}+D)\right)^{**}$ is a line bundle on $Y$ for every $m \geqq m_{0}$
(\cite[p.175,Proposition 4.2]{f-m}).  
\end{enumerate}
We note that adding effective exceptional $\mathbb{Q}$-divisors does not 
change the log canonical ring.   
Similarly as (\ref{HL}) we define the $\mathbb{Q}$-line bundle $L_{X/Y,D}$ on $Y$ by 
\begin{equation}\label{logHL}
L_{X/Y,D} = \frac{1}{m!}\,f_{*}\mathcal{O}_{X}(m!(K_{X/Y}+D))^{**}. 
\end{equation}
$L_{X/Y,D}$ is independent of the choice of a sufficiently large $m$ (\cite[p.169,Proposition 2.2]{f-m}). 
Let us fix such a $m$. 
Similarly as before we shall define the singular hermitian metric on $L_{X/Y,D}$ by 
\begin{equation}
h_{L_{X/Y,D}}^{m!}(\sigma,\sigma)(y):= \left(\int_{X_{y}}|\sigma|^{\frac{2}{m!}}\right)^{m!},
\end{equation}
where $y \in Y - \Sigma$, $X_{y}:= f^{-1}(y)$ and $\sigma\in m!L_{X/Y,D,y}$. 
We note that since $(X,D)$ is KLT, $h_{L_{X/Y,D}}$ is well defined.
As before $h_{L_{X/Y,D}}$ has semipositive curvature in the sense of current
(cf. \cite{subad,b-p}). 
By the same strategy as in the proof of Theorem \ref{main}, we have  
the following KLT version of Theorem \ref{canmeasure}.  
\begin{theorem}\label{logcanmeasure}(\cite{canonical,LC}) 
In the above notations, there exists a unique singular hermitian metric 
on $h_{K}$ on $K_{Y}+L_{X/Y,D}$ and a nonempty Zariski open subset 
$U$ of $Y$ such that 
\begin{enumerate}
\item[\em (1)] $h_{K}$ is an AZD of $K_{Y} + L_{X/Y,D}$, 
\item[\em (2)] $f^{*}h_{K}$ is an AZD of $K_{X}+D$,  
\item[\em (3)] $h_{K}$ is $C^{\infty}$ on $U$, 
\item[\em (4)] $\omega_{Y} = \sqrt{-1}\,\Theta_{h_{K}}$ is a K\"{a}hler form 
on $U$,  
\item[\em (5)] $-\mbox{\em Ric}_{\omega_{Y}} + \sqrt{-1}\,\Theta_{h_{L_{X/Y,D}}} = \omega_{Y}$
holds on $U$. \fbox{} 
\end{enumerate}
\end{theorem}
\begin{remark}
In Theorem \ref{logcanmeasure}, the metric $h_{K}$  depends only on 
the logcanonical ring of $(X,D)$.  Hence adding effective exceptional
$\mathbb{Q}$-divisors does not affect $h_{K}$ and $\omega_{Y}$ essentially. 
\fbox{} 
\end{remark}

\noindent We define {\bf the canonical measure $d\mu_{can}$ of the KLT pair $(X,D)$} by 
\begin{equation}
d\mu_{can} := \frac{1}{n!}f^{*}\left(\omega_{Y}^{n}\cdot h_{L_{X/Y,D}}^{-1}\right), 
\end{equation} 
where $n = \dim Y$. 
 
\subsection{Relative canonical measures}\label{relativemeasure}

In the previous subsection, we have introduced the (log) canonical measure
on a KLT pair $(X,D)$ with nonnegative Kodaira dimension. 
In this subsection, we consider the  variation of 
canonical measures on an algebraic fiber space.     
Let $f : X\to Y$ be an algebraic fiber space and let $D$ be an effective 
$\mathbb{Q}$-divisor such that $(X,D)$ is KLT.   Let $Y^{\circ}$ denote 
the complement of the discriminant locus of $f : X \to Y$. 
For a general $y\in Y^{\circ}$, $(X_{y},D_{y})$ is a KLT pair.
We denote the set : $\{y\in Y^{\circ}|(X_{y},D_{y})\,\,\mbox{is KLT}\}$ by $Y_{0}$.  We assume that 
$\mbox{Kod}(X_{y},D_{y}) \geqq 0$ holds for $y\in Y_{0}$.  By \cite{can}, we see that $h^{0}(X_{y},\mathcal{O}_{X_{y}}(m(K_{X_{y}}+D_{y})))$ is constant over 
$Y_{0}$ for every $m > 0$ such that $mD$ is Cartier and $f_{*}\mathcal{O}_{Y}(m(K_{X/Y}+D))$ is locally free over $Y_{0}$ for such a $m$. 
Then by  Theorem \ref{logcanmeasure},
we may define the canonical measure $d\mu_{can,y}$ of $(X_{y},D_{y})$. 
The family $\{ d\mu_{can,y}^{-1}\}_{y\in Y_{0}}$ defines a singular hermitian 
metric $h_{K}$ on $K_{X/Y} + D$.   The following theorem asserts that 
$h_{K}$ has semipositive curvature\footnote{Of course the main assertion is 
the semipositivity of the curvature in horizontal direction 
with respect to $f : X \to Y$}.

\begin{theorem}(\cite[Theorem 4.1]{canonical})\label{relative}
Let $f : X \longrightarrow Y$ be an algebraic fiber space.  And let $D$ be an effective divisor 
on $X$ such that $(X,D)$ is KLT.  
Suppose that $f_{*}\mathcal{O}_{Y}\left(\lfloor m(K_{X/Y}+D)\rfloor\right) \neq 0$ for some $m > 0$.
Then there exists a singular hermitian metric $h_{K}$ on $K_{X/Y}+D$ 
such that 
\begin{enumerate}
\item[(1)] $\omega_{X/Y}:= \sqrt{-1}\,\Theta_{h_{K}}$ is semipositive on $X$,
\item[(2)] For a general smooth fiber $X_{y} := f^{-1}(y)$ such that 
$(X_{y},D_{y})$ is KLT,  
$h_{K}|X_{y}$ is $d\mu_{can,(X_{y},D_{y})}^{-1}$, where $d\mu_{can,(X_{y},D_{y})}$ denotes 
the canonical measure on $(X_{y},D_{y})$.  
In particular $\omega_{X/Y}|X_{y}$ is the canonical semipositive current on $(X_{y},D_{y})$ 
constructed as in Theorem \ref{main}.\fbox{} 
\end{enumerate} 
\end{theorem}
We call   
\begin{equation}
d\mu_{can,(X,D)/Y} : = h_{K}^{-1}
\end{equation}
{\em the relative log canonical measure} for the family of KLT pairs $f : (X,D) \to Y$. 
Theorem \ref{relative} is the direct consequence of the dynamical construction 
of the canonical measures (cf. \cite[Theorem 1.7]{canonical}) and the 
plurisubharmonic variation property of Bergman kernels (\cite{b2,KE,b-p}).  

\subsection{Weak semistability}\label{ws}

In this subsection we  prove the 2nd and the 3rd assertions 
in Theorems \ref{main} and \ref{logmain}. 
The proof follows closely the one of Theorem \ref{viehweg} in \cite{v}. 
But we replace the use of branched coverings  in \cite{v} by the use of Theorem \ref{logcanmeasure}.  This enables us to get rid of the assumption that 
$K_{X/Y}$ is $f$-semiample.  

Let us start the proof.  
Let $f : X \longrightarrow Y$ be an algebraic fiber space. 
And let $Y^{\circ}$ be the complement of the discriminant locus of $f$.
And let $X^{\circ}:= f^{-1}(Y^{\circ})$. 
We set  $r = \mbox{rank}\,f_{*}\mathcal{O}_{X}(mK_{X/Y})$
and let $X^{r}:= X\times_{Y}X\times_{Y}\cdots\times_{Y}X$ denote the $r$-times fiber product over $Y$ and let $f^{r}: X^{r}\to Y$ be the natural morphism. 
Then we have the natural morpshim:  
\begin{equation}
\det f_{*}\mathcal{O}_{X}(mK_{X/Y}) \rightarrow 
\otimes^{r}f_{*}\mathcal{O}_{X}(mK_{X/Y}) = f^{r}_{*}\mathcal{O}_{X^{r}}(mK_{X^{r}/Y}).    
\end{equation}
Hence we have the canonical global section  
\begin{equation}
\gamma \in \Gamma \left(X,f^{r*}(\det f_{*}\mathcal{O}_{X}(mK_{X/Y}))^{-1}
\otimes \mathcal{O}_{X^{r}}(mK_{X/Y})\right). 
\end{equation}
Let $\Gamma$ denote the zero divisor of $\gamma$.  
It is clear the $\Gamma$ does not contain any fiber over $Y^{\circ}$. 
Now we set 
\begin{equation}\label{delta0}
\delta_{0}:= \sup \{\delta > 0| \mbox{$(X^{r}_{y},\delta\cdot \Gamma_{y})$
is KLT for every $y\in Y^{\circ}$}\}.
\end{equation}
Let us take a positive rational number $\varepsilon < \delta_{0}$. 
Then we have  that there exists the relative  canonical measure
$d\mu_{can,(X^{r},\Delta)}$ on 
$f: (X^{r},\Delta)\longrightarrow Y$ as in Theorem \ref{logcanmeasure}. 
By the logarithmic plurisubharmonicity  of the canonical measure (Theorem \ref{relative}),
we see that 
\begin{equation}
\sqrt{-1}\partial\bar{\partial}\log d\mu_{can,(X^{r},\Delta)/Y} \geqq 0
\end{equation}
holds on $X$ in the sense of current.
We set 
\begin{equation}\label{sect2}
H_{m,\varepsilon} := d\mu_{can,(X^{r},\Delta)/Y}^{-1}. 
\end{equation}
Then $H_{m,\varepsilon}$ is a singular hermitian metric on 
\begin{equation}
(1+m\varepsilon)K_{X^{r}/Y} - \varepsilon\cdot f^{r*}\det f_{*}\mathcal{O}_{X}(mK_{X/Y})
\end{equation}
with semipositive curvature current by Theorem \ref{relative} 
and $H_{m,\varepsilon}|X^{r}_{y}$ is an AZD of 
\begin{equation}
(1+m\varepsilon)K_{X^{r}_{y}} - \varepsilon\cdot f^{r*}\det f_{*}\mathcal{O}_{X}(mK_{X/Y})|X_{y}
\end{equation}
for every $y\in Y^{\circ}$. 
Hence by \cite{b-p}, we have that 
\begin{equation}
f^{r}_{*}\mathcal{O}_{X^{r}}(K_{X^{r}/Y}+\ell(1+m\varepsilon)K_{X^{r}/Y})
\succeq \ell\varepsilon \det f_{*}\mathcal{O}_{X}(mK_{X/Y})
\end{equation}
holds for every positive integer $\ell$ such that $\ell\varepsilon$ is 
an integer. 
Since
\begin{equation}
f^{r}_{*}\mathcal{O}_{X^{r}}(K_{X^{r}/Y}+\ell(1+m\varepsilon)K_{X^{r}/Y})
= f_{*}\mathcal{O}_{X}(r(1+\ell(1+m\varepsilon))K_{X/Y})
\end{equation}
holds, we have that 
\begin{equation}\label{prec}
 f_{*}\mathcal{O}_{X}(r(1+\ell(1+m\varepsilon))K_{X/Y})
 \succeq \ell\varepsilon \det f_{*}\mathcal{O}_{X}(mK_{X/Y}) 
\end{equation}
holds.  By \cite{b-c-h-m}, we have that for every sufficiently large 
integer $a$,  the natural morphism:
\begin{equation}\label{fg}
\otimes^{k}f_{*}\mathcal{O}_{X}(a!K_{X/Y}) \to f_{*}\mathcal{O}_{X}(ka!K_{X/Y})
\end{equation}
is surjective for every $k\geqq 0$. 
Hence dividing the both sides of (\ref{prec}) by  $\ell(1+m\varepsilon)$ and letting $\ell$ tend to infinity, by the surjection (\ref{fg}) we have that for every sufficiently large positive integer $d$, 
\begin{equation}
f_{*}\mathcal{O}_{X}(d!K_{X/Y})\succeq \frac{d!\,\varepsilon}{(1 + m\varepsilon)r}\det f_{*}\mathcal{O}_{X}(mK_{X/Y})^{**}
\end{equation}
holds. \fbox{}

\section{Moduli spaces of metrized canonical models}

So far we have completed the proof of  the 2nd and the 3rd assertions in Theorems \ref{main} and \ref{logmain} (cf. Section \ref{ws}).  To prove the 1st assertion of Theorems \ref{main} or \ref{logmain}, we need to use the moduli space of metrized canonical models. 

In this section, we shall construct the moduli space of metrized canonical models (cf. Definition \ref{metrizedcanonical}) and prove that it is an algebraic space in the sense of \cite{ar}.   
Here we shall explain only the abosolute case, i.e., we do not explain the 
case of KLT pairs for simplicity.  The general case follows from the similar 
argument.   Hence we  omit it. 

\subsection{Metrized canonical models}

Let $X$ be a smooth projective variety with $\mbox{Kod}(X)\geqq 0$. 
By \cite{b-c-h-m}, we see that the canonical ring:  
 $R(X,K_{X}) := \oplus_{m=0}^{\infty}\Gamma (X,\mathcal{O}_{X}(mK_{X}))$ 
is finitely generated.  
Then 
\begin{equation}
Y := \mbox{Proj}\,R(X,K_{X})
\end{equation}
is called the canonical model of $X$.  Then $Y$ has only canonical singularities  and the Hodge $\mathbb{Q}$-line bundle $L_{X/Y}$ is defined on $Y$ (cf. 
Section 2.3).        
Unless $X$ is of general type, the  canonical model $Y$ does not reflect the full informaion of the canonical ring $R(X,K_{X})$. The full information of 
the canonical ring is recovered from $K_{Y}$ and $L_{X/Y}$ by the isomorphism: 
\begin{equation}\label{caniso2}
R(X,K_{X}) \simeq R(Y,K_{Y}+L_{X/Y})^{(a)}, 
\end{equation}
where $a$ is the minimal positive integer such that $f_{*}\mathcal{O}_{X}(aK_{X/Y})\neq 0$.  by using the Hodge $\mathbb{Q}$-line bundle $L_{X/Y}$.  
Hence it is natural to consider the pair $(Y,L_{X/Y})$ instead of $Y$. 
But to describe the semipositivity of $f_{*}\mathcal{O}_{X}(mK_{X/Y})$, even the pair 
$(Y,L_{X/Y})$ is not enough.  
Hence we consider the triple $(Y,(L_{X/Y},h_{L_{X/Y}}))$ instead of  $Y$, where $(L_{X/Y},h_{L_{X/Y}})$ is the Hodge $\mathbb{Q}$-line bundle on $Y$ with the Hodge metric $h_{L_{X/Y}}$ (cf. Section \ref{hodgebundle}).  We call the pair $(L_{X/Y},h_{L_{X/Y}})$ the metrized Hodge $\mathbb{Q}$-line bundle. 
Then we may recover the canonical ring $R(X,K_{X})$ from $(Y,(L_{X/Y},h_{L_{X/Y}}))$  by the isomorphism (\ref{caniso2}).  
Moreover we may recover the canonical K\"{a}hler current $\omega_{Y}$ by
solving the equation: 
\begin{equation}
-\mbox{Ric}_{\omega_{Y}} + \sqrt{-1}\,\Theta_{h_{L_{X/Y}}} = \omega_{Y}
\end{equation}
in terms of the dynamical construction as in \cite{canonical} (cf. Theorem \ref{DS} below).
\begin{definition}\label{metrizedcanonical}
The pair $(Y,(L_{X/Y},h_{L_{X/Y}}))$ above is said to be the {\bf metrized canonical model}
of $X$. \fbox{}
\end{definition} 
Hereafter we shall construct the moduli space of the metrized canonical models. 

\subsection{Construction of  the moduli space and the statement of the result}
Let $f : X \to Y$ be an algebraic fiber space and let $Y^{\circ}$ denote 
the complement of the discriminant locus of $f$.  
\begin{equation*}
\begin{picture}(400,400)
\setsqparms[1`1`1`1;350`350]
\putVtriangle(0,00)%
[X` Z ` Y;g `f `h]
\end{picture}
\end{equation*}  
be the relative Iitaka fibration such that 
$Z$ is the relative canonical model on $Y^{\circ}$ and 
$g : X \to Z$ is a morphism. 
Let $(L_{X/Z},h_{L_{X/Z}})$ be the Hodge $\mathbb{Q}$-line bundle on $Z$.
We consider the set 
\begin{equation}
\mathcal{U}:= \{(Z_{y},(L_{X/Z},h_{L_{X/Z}})|Z_{y})| y\in Y^{\circ}\}. 
\end{equation}
Let $a$ be the minimal positive integer such that 
$aL_{X/Z}$ is Cartier. 
We define the equivalence relation $\sim$ on $\mathcal{U}$ by 
\begin{equation}\label{equiv}
(Z_{y},(L_{X/Z},h_{L_{X/Z}})|Z_{y})\sim (Z_{y^{\prime}},(L_{X/Z},h_{L_{X/Z}})|Z_{y^{\prime}}), 
\end{equation}
if and only if there exists a biholomorphism: 
$\varphi : Z_{y}\to Z_{y^{\prime}}$ and a bundle isomorphism: 
$\tilde{\varphi} : aL_{X/Z} |Z_{y}\to aL_{X/Z}|Z_{y^{\prime}}$ such that
the following commutative diagram :  
\begin{equation*}
\begin{picture}(500,500)
\setsqparms[1`1`1`1;600`400]
\putsquare(0,00)%
[aL_{X/Z}|Z_{y}` aL_{X/Z}|Z_{y^{\prime}} ` Z_{y}`Z_{y^{\prime}};
\tilde{\varphi}` ` `\varphi]
\end{picture}
\end{equation*}  
and 
\begin{equation}
\tilde{\varphi}^{*}(h_{L_{X/Z}}|Z_{y^{\prime}}) = h_{L_{X/Z}}|Z_{y}
\end{equation}
hold. 
Then we define the set $\mathcal{M}$ by 
\begin{equation}
\mathcal{M} := \mathcal{U}/\sim   
\end{equation}
and call it the moduli space of metrized canonical models 
associated with $f : X \to Y$. 

At this moment it is not clear the $\mathcal{M}$ has a complex structure.
In this section we start to prove the following theorem. 

\begin{theorem}\label{quasiprojective}
The moduli space of metrized canonical models  $\mathcal{M}$ 
(associated with $f : X \to Y$) has a structure of  quasiprojective variety.  \fbox{}
\end{theorem}

\subsection{Topological structure on $\mathcal{M}$}

To endow the topology and the complex structure on $\mathcal{M}$, first we  identify $\mathcal{M}$ with a quotient of certain subset of a Hilbert scheme.  

\begin{lemma}\label{matsusaka} There exists a positive integer $m_{0}$ such that 
for every $m \geqq m_{0}$ and $(Z_{y},(L_{X/Z},h_{L_{X/Z}})|Z_{y})\in \mathcal{U}$,
the complete linear system $|am(K_{Z_{y}}+L_{X/Z})|$ embeds 
$Z_{y}$ into a projective space $\mathbb{P}^{N(m)}$, where 
$N(m)$ is a positive integer independent of $y\in Y^{\circ}$. \fbox{}
\end{lemma}

Let $m_{0}$ be a positive integer as in Lemma \ref{matsusaka} and let $m$ be a positive integer greater than or equal to $m_{0}$. 
Let $(Z_{y},(L_{X/Z},h_{L_{X/Z}})|Z_{y})\in \mathcal{U}$ be an arbitrary point and 
let $\omega_{Z_{y}}$ denote the canonical K\"{a}hler current on 
$Z_{y}$ (cf. Definition \ref{L-K-E}). 
Let $\{ \sigma_{0},\cdots ,\sigma_{N(m)}\}$ be an orthonormal basis of 
$H^{0}(Z_{y},\mathcal{O}_{Z_{y}}(am(K_{Y}+L_{X/Z}|_{Z_{y}})))$ with respect to 
the inner product:  
\begin{equation}
(\sigma,\sigma^{\prime}):= \int_{Z_{y}}h_{L_{X/Z}}^{am}(\omega_{Z_{y}}^{n})^{-(am-1)}\cdot\sigma\cdot\overline{\sigma^{\prime}}. 
\end{equation}
Let $[\Phi_{m}(Z_{y})]$ denote the Hilbert point corresponding to the 
embedding:
\begin{equation}
\Phi_{m}(z) := [\sigma_{0}(z):\cdots :\sigma_{N(m)}] (z\in Z_{y}). 
\end{equation}
We consider the set 
\begin{equation}
\mathcal{U}_{m}:= \{ [\Phi_{m}(Z_{y})]|(Z_{y},(L_{X/Z},h_{L_{X/Z}})|Z_{y})\in \mathcal{U}\}, 
\end{equation}
where $\Phi_{m}(Z_{y})$ runs all the choice of  orthonormal basis 
$\{ \sigma_{0},\cdots ,\sigma_{N(m)}\}$.
We set 
\begin{equation}
\mathcal{U}_{\infty} := \prod_{m=m_{0}}^{\infty}\mathcal{U}_{m}
\end{equation}
and 
\begin{equation}
G_{\infty} := \prod_{m=m_{0}}^{\infty}PU(N(m)+1),
\end{equation}
where for a positive integer $k$, $PU(k+1)$ denotes the projective unitary group acting on $\mathbb{P}^{k}$.

\begin{lemma}(\cite{ti,ze})\label{kernel}
Let $(Z_{y},(L_{X/Z},h_{L_{X/Z}})|_{Z_{y}})\in \mathcal{U}$ be an arbitrary point and 
let $\{ \sigma_{0},\cdots ,\sigma_{N(m)}\}$ be the orthonormal basis 
of $H^{0}(Z_{y},\mathcal{O}_{Z_{y}}(am(K_{Y}+L_{X/Z})))$ as above.  Then
the Bergman kernel:  
\begin{equation}
K_{am}:= \sum_{i=0}^{N(m)}|\sigma_{i}|^{2}
\end{equation}
satisfies the identity:
\begin{equation}
(\omega_{Z_{y}}^{n})^{-1}\cdot h_{L_{X/Z}}|_{Z_{y}} := \lim_{m\to\infty} K_{am}^{-\frac{1}{am}}
\end{equation}
compact uniformly with respect to the $C^{\infty}$-topology  on the complement of the discriminant locus of $g|_{X_{y}}: X_{y} \to  Z_{y}$. 
\fbox{} \vspace{3mm} 
\end{lemma} 
By Lemma \ref{kernel}, we have the natural identification:
\begin{equation}
\mathcal{M} := \mathcal{U}_{\infty}/G_{\infty},
\end{equation} 
where $G_{\infty}$ acts on $\mathcal{U}_{\infty}$ in the natural manner. 
Hence  $\mathcal{M}$ has a natural topological 
space structure with respect to the quotient topology.

\subsection{Complex structure on $\mathcal{M}$}

Although $\mathcal{U}_{m}$ does not have a natural complex structure apriori,
we may endow a natural complex structure on $\mathcal{M}$ using the 
variation of Hodge structure and the logarithmic deformation.   

The reason is that the Hodge $\mathbb{Q}$-line bundle  $(L_{X/Z},h_{L_{X/Z}})$ is
nothing but the pull back of  the universal 
line bundle on the period domain by the (reduced) period map. 
But since the Hodge line bundle $L_{X/Z}$ is not a genuine line bundle,  
we need to take a cyclic covering to define the period map.  
This makes the proof a little bit more complicated. 

First we shall define the period map on a family of a metrized canonical model.  Let $f : X \to Y$ be an algebraic fiber space with $\mbox{Kod}(X/Y) 
\geqq 0$. 
Let $g : X \to Z$ be the relative Iitaka fibration with respect to 
$f : X \to Y$ as above and we set 
\begin{equation}
k:= \dim X/Z = \dim X - \dim Z. 
\end{equation}  
Let $F_{z}$ denote the fiber of 
$g : X \to Z$ over $z\in Z$.   Let $Z^{\circ}$ denote the 
complement of the discriminant locus of $g : X \to Z$. 
Then $\mbox{Kod}(F_{z}) = 0$ holds for every $z\in Z^{\circ}$. 
Let $a$ be a minimal positive integer such that $|aK_{F_{z}}|
\neq \emptyset$ for every $z\in Z^{\circ}$. 
Then for every  $z\in Z^{\circ}$ there exists a  nonzero element  $\eta_{z}\in \Gamma (F_{z},\mathcal{O}_{F_{z}}(aK_{F_{z}}))$ 
and let 
\begin{equation}
\mu_{z}: \tilde{F}_{z}\to F_{z}
\end{equation}
be the normalization of the cyclic cover which uniformize 
$\sqrt[a]{\eta_{z}}$.   Let us consider the family $\{\tilde{F}_{z}\}_{z\in Z^{\circ}}$.  This family is not  well defined  over 
$Z^{\circ}$, but it defines a family
\begin{equation}
\tilde{f} : \tilde{X}^{\circ}\to \tilde{Z}^{\circ}
\end{equation}
over the finite unramified covering 
\begin{equation}
\varpi : \tilde{Z}^{\circ}\to Z^{\circ} 
\end{equation}
corresponding to the monodromy representation of the fundamental group 
\begin{equation}
\pi_{1}(Z^{\circ}) \to \mathbb{Z}/a\mathbb{Z}.
\end{equation}
We take a $\mathbb{Z}/a\mathbb{Z}$ equivariant resolution $Z^{(a)}\to \tilde{Z}^{\circ}$ and let 
\begin{equation}  
g^{(a)} : X^{(a)} \to Z^{(a)}
\end{equation}
be the resulting family of the cyclic $a$-coverings.  
We set 
\begin{equation}\label{zar}
U:= \mbox{the complement of the discriminant locus of $g^{(a)}$}
\end{equation}
and 
let 
\begin{equation}
g_{U}: (g^{(a)})^{-1}(U) \to U
\end{equation}
be the restriction of $g^{(a)}$.  
Let $\mathbb{E}\to  U$ be the local system $R^{k}g_{U*}\mathbb{C}$ 
and let $\{ \mathbb{F}^{p}\}_{p=0}^{k}$ be the Hodge filtration of 
$\mathbb{E}$.    
Then we have the period map 
\begin{equation}
\Phi : U \to  \Gamma\,\backslash D
\end{equation}
associated with the variation of Hodge structures, where $D$ is the period domain and $\Gamma$ denotes the image of the monodromy representation
of $\pi_{1}(U)$ to $\mbox{Aut}(D)$.  
In this geometric variation of Hodge structures, it is known that 
$\Gamma$ acts on $D$ properly discontinuously (\cite{griff}). 
Hence $\Gamma\,\backslash D$ is a complex space. 
Let $\overline{U}$ be the completion of $U$ such that the boundary 
$B : = \overline{U} - U$ is a divisor with normal crossings.  
Then by \cite{deligne}, the quasi canonical extension $\overline{\mathbb{E}}$ 
of $\mathbb{E}\otimes \mathcal{O}_{U}$ exists, i.e., $\overline{\mathbb{E}}$ is a locally free 
sheaf with the Gauss-Manin connection:  
\begin{equation}
\nabla : \overline{\mathbb{E}} \to \Omega^{1}_{\overline{U}}(\log B)
\otimes \overline{\mathbb{E}}
\end{equation}   
such that the real part of the eigenvalues of the residues around
components of $B$  lie in $[0,1)$. 
Since we have assumed that $B$ is a divisor with normal crossings, 
the Hodge filtration $\{\mathbb{F}^{p}\}$ extends as a filtration
$\{\mathcal{F}^{p}\}$ of 
$\overline{\mathbb{E}}$ by subbundles. 
Then the metrized Hodge $\mathbb{Q}$-line bundle $(L_{X/Z},h_{L_{X/Z}})$ 
corresponds to the Hodge bundle $\mathbb{F}^{k}$ induced by the period map $\Phi$.  Moreover the metric $h_{L_{X/Z}}$ is induced by the Hodge metric 
on the universal Hodge bundle on the period domain $D$. 
Here the Hodge metric is induced from the Hodge bilinear form.

Let $Z^{0}$ be the maximal Zariski open subset of $Z$ such that 
$h_{L_{X/Z}}|_{Z}$ is locally bounded.  We note that $Z^{0}$ may be much larger 
than the complement of the discriminant locus of $g : X \to Z$.
We note that $h_{L_{X/Z}}|U$ extends smoothly across the component $B_{i}$ such that 
the Picard-Lefschetz transformations on $\mathbb{F}^{k}$ are of finite order around $B_{i}$ 
(\cite{griff,sch}).  Let $B^{0}$ be the union of the irreducible 
components of $B$ such that the Picard-Lefschetz transforms 
on $\mathbb{F}^{k}$  around the 
components are of infinite order.  
Let $\varpi : \overline{U} \to Z$ be the natural morphism.
Then $Z^{0} = \varpi(\overline{U} - B^{0})$ holds (cf.\cite{sch}). 
We set 
\begin{equation}
S  =  \varpi (B^{0}). 
\end{equation}    
We consider the pair of the pairs:
\begin{equation}\label{quad}
\left((Z,S),\varpi_{*}(\overline{\mathbb{E}},\mathcal{F}^{k})\right). 
\end{equation}
Then by the above construction we have the following lemma.

\begin{lemma}\label{bijection} $\mathcal{U}:= \{ (Z_{y},(L_{X/Z},h_{L_{X/Z}})|Z_{y})| y\in Y^{\circ} \}$ is bijective to the set of quadruples: 
\begin{equation}
\mathcal{U}^{*}:= \{\left((Z_{y},S_{y}),\varpi_{y*}(\overline{\mathbb{E}}_{y},\mathcal{F}_{y}^{k})\right)| y\in Y^{\circ}\}.  
\end{equation}
\fbox{}
\end{lemma}
{\em Proof}.  Since $S_{y}$ is the polar locus of $h_{L_{X/Z}}|_{Z_{y}}$, 
$(Z_{y},(L_{X/Z},h_{L_{X/Z}})|Z_{y})$ determines  the pair $(Z_{y},S_{y})$. 
Since $(L_{X/Z},h_{L_{X/Z}})$ is determined by the period map:  
\[
\Phi : U \to \Gamma\,\backslash D,  
\]
$(Z_{y},(L_{X/Z},h_{L_{X/Z}})|Z_{y})$ determines the quadruple 
$\left((Z_{y},S_{y}),\varpi_{y*}(\overline{\mathbb{E}}_{y},\mathcal{F}_{y}^{k})\right)$. 
Conversely, since $\mathbb{E}_{y}$ is a flat vector bundle with the natural bilinear form,
the quadruple $\left((Z_{y},S_{y}),\varpi_{y*}(\overline{\mathbb{E}}_{y},\mathcal{F}_{y}^{k})\right)$ determines the pair $(Z_{y},(L_{X/Z},h_{L_{X/Z}})|Z_{y})$. 
This completes the proof. 
\fbox{} \vspace{3mm} \\
We define the equivalence relation $\sim$ on $\mathcal{U}^{*}$ by 
\[
\left((Z_{y},S_{y}),\varpi_{y*}(\overline{\mathbb{E}}_{y},\mathcal{F}_{y}^{k})\right) \sim 
\left((Z_{y^{\prime}},S_{y^{\prime}}),\varpi_{y^{\prime}*}(\overline{\mathbb{E}}_{y^{\prime}},\mathcal{F}_{y^{\prime}}^{k})\right)
\]
if and only if there exist a biholomorphism 
\begin{equation}
\varphi : (Z_{y},S_{y})\to 
(Z_{y^{\prime}},S_{y^{\prime}})
\end{equation}
and a sheaf isomorphism 
\begin{equation}
\tilde{\varphi} : (\varpi_{y})_{*}(\overline{\mathbb{E}}_{y},\mathcal{F}_{y}^{k})
\to (\varpi_{y^{\prime}})_{*}(\overline{\mathbb{E}}_{y^{\prime}},\mathcal{F}_{y^{\prime}}^{k}) 
\end{equation}
which covers $\varphi$ which induced by an isomorphism of the flat
vector bundles
\begin{equation}
\mathbb{E}|_{W_{y}}\to \mathbb{E}|_{W_{y^{\prime}}}, 
\end{equation}
where $W_{y},W_{y^{\prime}}$ are some nonempty Zariski open subsets 
of $U_{y}$ and $U_{y^{\prime}}$ (cf. (\ref{zar})) respectively. 

\begin{lemma}\label{cpx}
$\mathcal{M}^{*}:= \mathcal{U}^{*}/\sim$ has a structure of an algebraic space in the sense of \cite{ar}. 
\fbox{} 
\end{lemma}
{\em Proof of Lemma \ref{cpx}}.
Let $m_{0}$ be a sufficiently large positive integer such that 
$m_{0}!(K_{Z_{y}} + L_{X/Z}|_{Z_{y}})$ is Cartier and $|m_{0}!(K_{Z_{y}} + L_{X/Z}|_{Z_{y}})|$ 
is very ample for every $y\in Y^{\circ}$. 
We set $N := \dim  |m_{0}!(K_{Z_{y}} + L_{X/Z}|_{Z_{y}})|$.
If we  fix a basis 
of \\ $H^{0}(Z_{y},\mathcal{O}_{Z_{y}}(m_{0}!(K_{Z}+L_{X/Z}|_{Z_{y}})))$, 
then the basis gives an embedding: 
\[
\phi : Z_{y} \to \mathbb{P}^{N}. 
\]
and the images $\phi(Z_{y}),\phi(S_{y})$ define 
points in the Hilbert scheme $\mbox{Hilb}_{\mathbb{P}^{N}}$ of $\mathbb{P}^{N}$.  Hence  the linear system $|m_{0}!(K_{Z_{y}} + L_{X/Z}|_{Z_{y}})|$ gives 
an $PGL(N+1,\mathbb{C})$ orbit in $\mbox{Hilb}_{\mathbb{P}^{N}}\times 
\mbox{Hilb}_{\mathbb{P}^{N}}$. 
We denote the union of the orbits in $\mbox{Hilb}_{\mathbb{P}^{N}}\times 
\mbox{Hilb}_{\mathbb{P}^{N}}$ by $\mathcal{V}$ and let
\begin{equation}
\pi : (\mathcal{Z},\mathcal{S}) \to \mathcal{V}
\end{equation}   
be the universal family. 

Next we consider the pair $(\overline{\mathbb{E}}_{y},\mathcal{F}_{y}^{k})$.
Let $\mathcal{O}_{Z_{y}}(1)$ denote $\mathcal{O}_{Z_{y}}(m_{0}!(K_{Z_{y}}+L_{X/Z}|_{Z_{y}}))$.  For a positive integer $\ell$, we set  
\[
(\varpi_{y})_{*}\overline{\mathbb{E}}_{y}(\ell):= (\varpi_{y})_{*}\overline{\mathbb{E}}_{y}\otimes
\mathcal{O}_{Z_{y}}(\ell)\,\,\,\,\mbox{and}\,\,\,\, (\varpi_{y})_{*}\mathcal{F}_{y}^{k}(\ell):= (\varpi_{y})_{*}\mathcal{F}_{y}^{k}\otimes\mathcal{O}_{Z_{y}}(\ell). 
\]
Then for every $\ell \geqq 0$, we have the canonical inclusion:
\begin{equation}\label{subsp}
H^{0}(Z_{y}, (\varpi_{y})_{*}\mathcal{F}_{y}^{k}(\ell))
\hookrightarrow H^{0}(Z_{y},(\varpi_{y})_{*}\overline{\mathbb{E}}_{y}(\ell)).
\end{equation}
We denote $PGL(N+1,\mathbb{C})$ by $G$. 
If $(Z_{v},S_{v}),(Z_{v^{\prime}},S_{v^{\prime}}) (v,v^{\prime}\in\mathcal{V})$ are in the same orbit of $G$,  
then an element $g\in G$ induces a biholomorphism between 
$(Z_{v},S_{v})$ and $(Z_{v^{\prime}},S_{v^{\prime}})$
and an isomorphism of the flat vector bundles
$\mathbb{E}_{v}$ and $\mathbb{E}_{v^{\prime}}$ on the cyclic covers
and induces the isomorphism
between  $(\varpi_{v})_{*}\overline{\mathbb{E}}_{v}$ and $(\varpi_{v^{\prime}})_{*}
\overline{\mathbb{E}}_{v^{\prime}}$.  
The latter isomorphisms are unique up to the action of $\mathbb{Z}/a\mathbb{Z}$
and the $\mathbb{C}^{*}$-action. But since $\mathcal{F}^{k}_{v}(v\in \mathcal{V})$ is $\mathbb{Z}/a\mathbb{Z}$-equivariant subsheaf of $\mathbb{E}_{v}$, 
in spite of the umbiguity of the isomorphism, any such isomorphism maps the subspace 
$H^{0}(Z_{v}, (\varpi_{v})_{*}\mathcal{F}_{v}^{k}(\ell))
\subset H^{0}(Z_{v},(\varpi_{v})_{*}\overline{\mathbb{E}}_{v}(\ell))$
to the same subspace of 
$H^{0}(Z_{v^{\prime}}, (\varpi_{v^{\prime}})_{*}\mathcal{F}_{v^{\prime}}^{k}(\ell))$.  
We set 
\begin{equation}
\mathcal{M}^{\prime}:= \mathcal{V}/G. 
\end{equation}
Then by the construction, $\mathcal{M}^{\prime}$ is an algebraic space. 
And $\{(\varpi_{v})_{*}\overline{\mathbb{E}}_{v}(\ell)| v\in \mathcal{V}\}$ decends to a coherent sheaf $\mathcal{E}$ on $\mathcal{M}^{\prime}$ and 
the image of the inclusion (\ref{subsp}) determines the 
subsheaf $(\varpi_{v})_{*}\mathcal{F}^{k}$ for every $v\in \mathcal{V}$, if we take $\ell$ sufficiently large.  Let us fix such $\ell$.    
Hence by the properness of the period map (\cite{griff}), $\mathcal{M}^{*}$ is a locally closed subset (in Zariski topology) of the Grassmann bundle 
\begin{equation}
\mathcal{G}\to \mathcal{M}^{\prime} 
\end{equation}
associated with $\mathcal{E}$ with fiber $Gr(e,a)$, where 
$e = \mbox{rank}\,\mathcal{E}$. 
Hence $\mathcal{M}^{*}$ is an algebraic space. \fbox{} \vspace{3mm} \\

By Lemma \ref{bijection}, we see that there exists 
a homeomorphism between $\mathcal{M}$ and $\mathcal{M}^{*}$. 
Hence we have the complex structure on $\mathcal{M}$ by 
Lemma \ref{cpx}.

\subsection{Separatedness}

To ensure the existence of ${\cal M}$  as a Hausdorff complex space,
the following lemma is essential.
 
\begin{lemma}\label{separability}
Let $f : (X,D) \longrightarrow \Delta$ and $f : (X^{\prime},D^{\prime})\longrightarrow \Delta$
be flat projective families of KLT pairs with nonnegative Kodaira dimension 
over the unit open disk $\Delta$ in $\mathbb{C}$.
Let $\Delta^{*}:= \Delta -\{ 0\}$ denote the punctured disk.    And let $h : (Y,(L_{X/Y},h_{L_{X/Y}}))\to \Delta$, 
$h^{\prime} : (Y^{\prime},(L_{X/Y}^{\prime},h_{L_{X/Y}^{\prime}}))\to \Delta$ be 
the corresponding family of metrized pairs. 
Suppose that there exists an equivalence 
\begin{equation}
\varphi : (Y,(L_{X/Y},h_{L_{X/Y}}))|\Delta^{*}
\to (Y^{\prime},(L_{X/Y}^{\prime},h_{L_{X/Y}^{\prime}}))|\Delta^{*}
\end{equation}
of the families over $\Delta^{*}$ in the sense of (\ref{equiv}).
Then $\varphi$ extends uniquely to an equivalence between 
$(Y,(L_{X/Y},h_{L_{X/Y}}))$ and $(Y^{\prime},(L_{X/Y}^{\prime},h_{L_{X/Y}^{\prime}}))$.
\fbox{}
\end{lemma} 
{\em Proof of Lemma \ref{separability}}. 
Let $\varphi_{s}$ denote the restriction of $\varphi$ to $Y_{s}$ $(s\in \Delta^{*})$. 
Let $\omega_{s}$ denotes the canonical K\"{a}hler current on $Y_{s}$ 
constructed as in Theorem \ref{canmeasure}. 
Then  by the equation (\ref{LKE}), 
we see that $\varphi_{s} : Y_{s} \longrightarrow Y^{\prime}_{s} (s\in \Delta^{*})$ is an isometry between the K\"{a}hler spaces 
$(Y_{s},\omega_{Y_{s}})$ and $(Y^{\prime}_{s},\omega^{\prime}_{s})$.
Then by Ascoli-Arzela's theorem and Montel's theorem, we can easily see that 
$\varphi_{s}$ converges to an isometry 
\begin{equation}\label{isometry}
\varphi_{0} : (Y_{0,reg},\omega_{Y_{0}})\to  (Y_{0,reg}^{\prime},\omega_{Y_{0}^{\prime}})
\end{equation}
and  is holomorphic. This means that $\varphi$ extends uniquely to a biholomorphism between $Y$ and $Y^{\prime}$.  

The correspondence of the Hodge line bundles 
is obtained as follows.   By (\ref{isometry}), we have the equality:
\[
\varphi_{0}^{*}\,\,\omega_{Y_{0}^{\prime}}^{n} = \omega_{Y_{0}}^{n}. 
\]
Then by the equation (\ref{LKE}), we obtain that
\[
\varphi_{0}^{*}\,\,\Theta_{h_{L_{X/Y}^{\prime}}} = \Theta_{h_{L_{X/Y}}}
\]
holds on $Y_{0}$.  Hence we see that  
$\varphi$ extends uniquely to an equivalence between 
$(Y,(L_{X/Y},h_{L_{X/Y}}))$ and $(Y^{\prime},(L_{X/Y}^{\prime},h_{L_{X/Y}^{\prime}}))$.
\fbox{} \vspace{3mm} \\ 

By Lemma \ref{separability}, we see that $\mathcal{M}$ is separable.  
Then by the construction, we see that $\mathcal{M}$ is an separable 
algebraic space in the sense of Artin (cf. \cite{ar}).  
So far we have proven the following:
\begin{proposition}\label{algebraic}
$\mathcal{M}$ is a separable algebraic space. \fbox{}
\end{proposition}

\section{Descent of the Monge-Amp\`{e}re foliation}

Let $f : X \to Y$ be an algebraic fiber space such that $\mbox{Kod}(X/Y) \geqq 0$.   Then we have the relative canonical measure $d\mu_{can,X/Y}$ 
as in Theorem \ref{relative}.   Then by Theorem \ref{relative}, 
$\omega_{X/Y}:= \sqrt{-1}\,\Theta_{d\mu_{can,X/Y}^{-1}}$  is a closed 
semipositive current on $X$ which is generically $C^{\infty}$ 
by Theorem \ref{regularity} below.   
Then $\omega_{X/Y}$ defines a (possibly) singular foliation on  
$X$ whose leaves are complex analytic.  
In this section, we analyse this foliation.  

\subsection{Weak semistability and Monge Amp\`{e}re foliations}
Let $f : X \longrightarrow Y$ be a surjective projective morphism 
of smooth projective varieties with connected fibers such that 
$\mbox{Kod}(X/Y) \geqq 0$. 
Let $m$ be a  positive integer and let 
\begin{equation}
E_{m}:= f_{*}\mathcal{O}_{X}(mK_{X/Y}).
\end{equation}
We  assume that $E_{m}\neq 0$.  Let 
\begin{equation}
r:= \mbox{rank}\,\, E_{m}. 
\end{equation}
Let $h_{m}$ be 
the (singular) hermitian metric on $E_{m}$ defined by 
\begin{equation}\label{l2metric}
h_{m}(\sigma,\tau):= \int_{X/Y}\sigma\cdot\bar{\tau}\cdot h_{K,X/Y}^{m-1}. 
\end{equation}
Then since $h_{K,X/Y}|X_{y}$ is an AZD of $K_{X_{y}}$ for every $y\in Y^{\circ}$. We see that $h_{m}$ is a locally bounded hermitian metric on 
$E_{m}|Y^{\circ}$.    
$h_{m}$ defines an hermitian metric $\det h_{m}$ on $\det E_{m}$ 
and is locally bounded on $Y^{\circ}$. 
By \cite{KE} or \cite{b-p}, we see that $\sqrt{-1}\,\Theta_{\det h_{m}}$ is a closed positive current on $Y^{\circ}$.  

Let $X^{r}$ denote the $r$-times fiber product of $X$ over $Y$ 
and let 
\begin{equation}
f^{r}: X^{r} \to Y
\end{equation}
be the natural morphism. 
Let $\delta_{0}$ be the positive number as in Section 2 (cf. (\ref{delta0})) and let  $\varepsilon$ 
be a positive rational number such that $\varepsilon < \delta_{0}$.
Let $H_{m,\varepsilon}$ be the singular hermitian metric on
\begin{equation}
(1 + m\varepsilon)K_{X^{r}/Y} - \varepsilon f^{r*}\det h_{m}
\end{equation}
constructed as in Section 2 (cf. (\ref{sect2})). 
We define the singular hermitian metric  $H^{+}_{m,\varepsilon}$ 
on $K_{X^{r}/Y}$ by 
\begin{equation}
H^{+}_{m,\varepsilon}:= \left(H_{m,\varepsilon}\cdot f^{r*}(\det h_{m})^{\varepsilon}\right)^{\frac{1}{1+m\varepsilon}}. 
\end{equation}  
Since 
\begin{equation}
\Theta_{H^{+}_{m,\varepsilon}} = \frac{1}{1+m\varepsilon}
\left(\Theta_{H_{m,\varepsilon}} + \varepsilon\cdot f^{r*}\Theta_{\det h_{m}}\right)
\end{equation}
and $\sqrt{-1}\,\Theta_{\det h_{m}}$ is semipositive current on $Y$, 
we have the following lemma.
\begin{lemma}\label{incl}
\begin{equation}
\sqrt{-1}\Theta_{H^{+}_{m,\varepsilon}}
\geqq \frac{\varepsilon}{1 + m\varepsilon}\sqrt{-1}f^{r*}\Theta_{h_{\det h_{m}}}
\end{equation}
holds on $X$. \fbox{}
\end{lemma}
We set 
\begin{equation}
\omega_{m,\varepsilon} := \sqrt{-1}\,\Theta_{H_{m,\varepsilon}^{+}}.
\end{equation} 
Let $d\mu_{can,X^{r}/Y}$ be the relative canonical measure 
on the algebraic fiber space $f^{r} : X^{r}\to Y$.
We set 
\begin{equation}
\omega_{X^{r}/Y}:= \sqrt{-1}\,\partial\bar{\partial}\log\, d\mu_{can,X^{r}/Y}
\end{equation}  
$d\mu_{can,X^{r}/Y}$ is $C^{\infty}$ on a nonempty Zariski open 
subset $U$ of $X^{r}$ by Theorem \ref{regularity} below.   
Then we see that 
\begin{equation}
\mathcal{F}:= 
\{\xi \in TX^{r}|U\,\,; \,\, \omega_{X^{r}/Y}(\xi,\bar{\xi}) = 0\}
\end{equation}
defines a singular  foliation on an open subset $V$ of $U$ defined by 
\begin{equation}\label{vv}
U_{0}: = \{ x\in U|\mbox{$\omega_{X^{r}/Y}$ is of maximal rank at $x$}\},
\end{equation}
i.e., $\mathcal{F}$ is a Monge-Amp\`{e}re foliation 
associated with the semi K\"{a}hler form $\omega_{X^{r}/Y}|U_{0}$.
Hence $\mathcal{F}$ has  complex analytic leaves on $U_{0}$ (\cite{b-k}).  
But at this moment it is not clear $\mathcal{F}|U_{0}$ is a complex analytic 
foliation. 
By Lemma \ref{incl}, we have the following lemma.
 
\begin{lemma}\label{vanishing}
$f^{r*}\Theta_{\det h_{m}}|\mathcal{F} \equiv 0$ holds. \fbox{}
\end{lemma}
{\em Proof}. 
Since $d\mu_{can,X^{r}/Y}^{-1}$ is an AZD of $K_{X^{r}/Y}$,  we see that 
for every leaf $F$ of $\mathcal{F}$, 
\begin{equation}\label{vanish}
\Theta_{H_{m,\varepsilon}}|F \equiv 0
\end{equation}
holds.   

In fact otherwise,  we have a singular hermitian metric: 
\begin{equation}
H_{m,\varepsilon}^{1/2}\cdot d\mu_{can,X^{r}/Y}^{-1/2}
\end{equation}
on $K_{X^{r}/Y}$ with semipositive curvature and strictly bigger numerical dimension
than $d\mu_{can,X^{r}/Y}^{-1}$ . 
This contradicts  the fact that $d\mu_{can,X^{r}/Y}^{-1}$ is an AZD of 
$K_{X^{r}/Y}$.   

Hence combining (\ref{vanish}) and Lemma \ref{incl}, we see that 
\begin{equation}
f^{r*}\Theta_{\det h_{m}}|F \equiv 0
\end{equation}
holds.  This completes the proof of Lemma \ref{vanishing}.  \fbox{}

\subsection{Trivialization along the leaves on $Y$}
Let $(E_{m},h_{m})$ be as above.   Then  for any local
holomorphic section  $\xi$ of $E_{m}$ on some open subset $V$ of $Y$,
\begin{equation}
h_{m}(\sqrt{-1}\Theta_{h_{m}}(\xi),\xi) 
\end{equation}
is a semipositive $(1,1)$-current on $V$  by \cite{KE} or \cite{b-p}. 
Then the curvature $\det E_{m}$ is computed as:
\begin{equation}
\Theta_{\det h_{m}}(y) = \sum_{\alpha} h_{m}(\Theta_{h_{m}}(\mbox{\bf e}_{\alpha}),\mbox{\bf e}_{\alpha}), 
\end{equation}
where $\{\mbox{\bf e}_{\alpha}\}$ is an orthonormal basis of $E_{m,y}$
with respect to $h_{m}$. 
Hence $\sqrt{-1}\Theta_{\det h_{m}}$ is a closed semipositive current on $Y$.  
Since 
\begin{equation}
f^{r*}\Theta_{\det h_{m}}|\mathcal{F} \equiv 0
\end{equation}
holds by Lemma \ref{vanishing}, $\left(f^{r*}E_{m},f^{r*}h_{m}\right)$ is flat along 
every leaf of $\mathcal{F}$.   
Hence this implies that for every $x \in U_{0}$ and an orthonormal basis 
$\{\mbox{\bf e}_{\alpha,x}\}$ of $\left(f^{r*}E_{m}\right)_{x}$
with respect to $f^{r*}h_{m}$, the parallel transport of 
$\{\mbox{\bf e}_{\alpha,x}\}$  along the leaf $F$ of $\mathcal{F}$ containing $x$
 trivialize $(f^{r*}E_{m})|F$  locally. 
Let 
\begin{equation*}
\begin{picture}(400,400)
\setsqparms[1`1`1`1;350`350]
\putVtriangle(0,00)%
[X` Z` Y;g `f `h]
\end{picture}
\end{equation*}  
be the relative Iitaka fibration such that $Z$ is the family of relative 
canonical models and let $(L,h_{L})$ be the Hodge $\mathbb{Q}$-line bundle on $Z$.  
  
Then we have the following lemma : 
\begin{lemma}\label{trivial}
For every leaf $F$ of $\mathcal{F}$, the restriction 
\begin{equation}
(Z,(L,h_{L}))|f^{r}(F) \to f^{r}(F)
\end{equation}
is locally trivial. \fbox{}
\end{lemma} 
{\em Proof of Lemma \ref{trivial}}. 
By the flatness of $(E_{m},h_{m})$ along $f^{r}(F)$, we see that 
the parallel transport in $(E_{m},h_{m})|f^{r}(F)$ locally trivialize 
$E_{m}$ as above.  This implies that $Z|f^{r}(F)$ is also trivialized by the 
parallel transport, since it is the (log) canonical image.  Hence $(L,h_{L})|f^{r}(F)$ is  also locally trivial (as a metrized family of complex lines).  \fbox{}  

\subsection{Closedness of leaves}
Let $\mathcal{M}$ be the moduli space which parametrizes
the equivalence classes of   
\begin{equation}
\{(Z_{y},(L,h_{L})|Z_{y})| y \in Y^{\circ}\}
\end{equation}  
constructed as in Section 3. 
Now we consider the moduli map 
\begin{equation}\label{modmap}
\mu : Y^{\circ} \to \mathcal{M}
\end{equation}
defined by 
\begin{equation}
\mu(y):= [(Z_{y},(L,h_{L})|Z_{y}))], 
\end{equation}
where $[(Z_{y},(L,h_{L})|Z_{y}))]$ denotes the equivalence class in $\mathcal{M}$. 
Then by Lemma \ref{trivial} for every leaf $F$ of $\mathcal{F}$, 
$f^{r}(F)$ is contained in the fiber of $\mu : Y ^{\circ}\to \mathcal{M}$.
But by the construction, conversely, we see that for every 
$P\in \mathcal{M}$, $(f^{r})^{-1}(\mu^{-1}(P))$ is contained in a 
leaf of $\mathcal{F}$. 

Hence we conclude that for every leaf $F$ of $\mathcal{F}$, 
$f^{r}(F)$ is an open subset of the fiber of $\mu$ and 
$f_{*}\mathcal{F}$ descends to the foliation defined by 
the moduli map $\mu$.    Hence we may take $U_{0}$ defined as (\ref{vv}) to be a nonempty Zariski open subset of $X$.  
By the above argument we have the following lemma.
 
\begin{lemma}\label{closedness}
In the above notations, we have the followings:
\begin{enumerate}
\item[(1)] $\mathcal{F}$ decends to the foliation  $df(\mathcal{F})$ on $Y^{\circ}$, 
\item[(2)] Every leaf of $df(\mathcal{F})$ is closed in $Y^{\circ}$ 
and is a fiber of the moduli map $\mu : Y^{\circ}\to \mathcal{M}$,
\item[(3)] $\mathcal{F}$ is a singular analytic foliation on $X$. 
\fbox{}
\end{enumerate}
\end{lemma} 
{\em Proof.}  The assertions (1) and (2) have already been proven. 
The assertion (3) follows from (2) and Lemma \ref{trivial}.  \fbox{}

\section{Completion of the proof of Theorems \ref{main},\ref{logmain}
and \ref{quasiprojective}}

In this section we complete the proof of the proof of Theorems \ref{main},\ref{logmain} and \ref{quasiprojective}. 
But we shall omit the proof of Theorem \ref{logmain}, since the proof is essentially the same as the one of Theorem \ref{main}.   

Let $f : X \longrightarrow Y$ be an algebraic fiber space. 
Suppose that for a general fiber $F$ of $f$, $\mbox{Kod}(F) \geqq 0$
holds.  
Then we have the relative Iitaka fibration:  
\begin{equation*}
\begin{picture}(400,400)
\setsqparms[1`1`1`1;350`350]
\putVtriangle(0,00)%
[X` Z ` Y;g `f `h]
\end{picture}
\end{equation*}  
such that $Z$ is a family of relative canonical models. 
By taking a suitable modification of $X$, we may assume that 
$g$ is a morphism. 

Let $(L,h_{L})$ be the Hodge line bundle on $Z$ as in Section \ref{hodgebundle}. 
Then we have that 
\begin{equation}
f_{*}\mathcal{O}_{X}(m!K_{X/Y}) = h_{*}\mathcal{O}_{Y}(m!(K_{Z/Y} + L))
\end{equation}
holds for every sufficiently large $m$.  
Let $Y^{\circ}$ be the complement of the discriminant locus of 
$f: X\to Y$. 
Let 
\begin{equation}
\mu : Y^{\circ} \to \mathcal{M}.
\end{equation}
be the moduli map (\ref{modmap}) as above. 
Then by the quasi-unipotence of the monodromy (\cite{borel}), 
we see that there exists a positive integer $b$ such that
for every $m > 0$  
\begin{equation}
\left(\det f_{*}\mathcal{O}_{X}(mK_{X/Y})\right)^{\otimes b}
\end{equation}
and 
\begin{equation}
\left(f_{*}\mathcal{O}_{X}(mK_{X/Y})\right)^{\otimes b}
\end{equation}
decend to vector bundles on $\mathcal{M}$.
Then the relative canonical measure $d\mu_{can,X/Y}$ defines a $L^{2}$ 
metric $h_{m}$ on $f_{*}\mathcal{O}_{X}(mK_{X/Y})$ as in (\ref{l2metric}) and then 
$h_{m}$ defines a singular hermitian metric $\det h_{m}$ on 
$\det f_{*}\mathcal{O}_{X}(mK_{X/Y})$.  The metric $h_{m}$ is an invariant 
metric by Theorem \ref{canmeasure}.  In the above notations, we have 
the following lemma. 
 
\begin{lemma}\label{numpos}
Let $a$ be the minimal positive integer such that 
$f_{*}\mathcal{O}_{X}(aK_{X/Y}) \neq 0$. 
Let $m_{0}$ be a sufficiently large positive integer.
Then 
\begin{equation}
F:= \mu_{*}\left(\det f_{*}\mathcal{O}_{X}(m_{0}aK_{X/Y})\right)^{\otimes b}
\end{equation}
is a line bundle on $\mathcal{M}$ with the hermitian metric $h_{F}$ such that
\begin{enumerate} 
\item[\em (1)] $\mu^{*}h_{F} = h_{am_{0}}$,
\item[\em (2)] For every subvariety $V$ in $\mathcal{M}$, 
$(F|_{V},h_{F}|_{V})$ is big on $V$ (cf. Definition \ref{numerical}).  
\end{enumerate}
\fbox{}  
\end{lemma}
{\em Proof of Lemma \ref{numpos}}. 
The first  assertion (1) is trivial by the construction and the birational 
invariance of the canonical measures. 

By Lemma \ref{closedness} we have the followings:
\begin{enumerate}
\item[(1)] The foliation $\mathcal{F}$ decends to a foliation $df(\mathcal{F})$ on $Y$.
\item[(2)] $\omega_{Z/Y}$ is generically strictly positive in the 
 transverse direction with respect to $\mathcal{F}$.
\item[(3)] $\mu$ contracts the leaf of $df(\mathcal{F})$. 
\end{enumerate}
Then the second assertion (2) holds, if $V = \mathcal{M}$ by the construction, 
  
For a general $V$, the assertion (2) follows from the functoriality.  
\fbox{}\vspace{3mm} \\

\noindent By Proposition \ref{algebraic}, we see that $\mathcal{M}$ has 
a structure of a separable algebraic space. 
Then by Lemma \ref{numpos} and the quasiprojectivity criterion Theorem \ref{criterion} below, we see that $\mathcal{M}$ is quasiprojective. 
This completes the proof of Theorem \ref{quasiprojective}.  

To complete the proof of Theorem \ref{main} we use the 
weak semipositivity (cf. (\ref{wp1}) or (3)(b)) in Theorem \ref{main}. Then we see that 
$\mu_{*}\left(f_{*}\mathcal{O}_{X}(mK_{X/Y})\right)^{\otimes b}$
is globally generated on $\mathcal{M}$ for every sufficiently large 
and divisible $m$.  
Then since  
\begin{equation}
\mu^{*}\left(\mu_{*}\left(f_{*}\mathcal{O}_{X}(mK_{X/Y})\right)^{\otimes b}\right) = \left(f_{*}\mathcal{O}_{X}(mK_{X/Y})\right)^{\otimes b}
\end{equation} 
holds by the construcion, we see that  $\left(f_{*}\mathcal{O}_{X}(mK_{X/Y})\right)^{\otimes b}$
is globally generated on $Y^{\circ}$ for every sufficiently large $m$. 
Then by the finite generation of canonical rings (\cite{b-c-h-m}), this implies that 
there exists a positive integer $m_{0}$ such that 
$f_{*}(mK_{X/Y})$ is globally generated over 
$Y^{\circ}$, if $b|m$ and $m \geqq m_{0}$.  
This completes the proof of Theorem \ref{main}.  
The proof of Theorem \ref{logmain} is similar.  \fbox{}

\section{Parameter dependence of canonical measures} 

In this section we  prove the following regularity theorem 
for the relative canonical measure $d\mu_{can,X/Y}$ constructed as in 
Theorem \ref{relative}.

\begin{theorem}\label{regularity}
Let $f : X \to Y$ be an algebraic fiber space with $\mbox{\em Kod}(X/Y)\geqq 0$. Let $d\mu_{can,X/Y}$ be the relative canonical measure on 
$X$ constructed as in Theorem \ref{relative}. 
Then $d\mu_{can,X/Y}$ is $C^{\infty}$ on a nonempty Zariski open subset 
of $X$. \fbox{}
\end{theorem}
\begin{remark}
Essentially the same  regularity result holds for the relative log canonical measure 
$d\mu_{can,(X,D)/Y}$ for the family of KLT pairs $f : (X,D)\to Y$ (cf. Theorem \ref{relative}).   The proof requires the dynamical construction of 
(log) canonical measures as in \cite{LC}, but otherwise the proof is the same as the one of Theorem \ref{regularity}.  \fbox{}
\end{remark}
Here I would like to explain the scheme of the proof of Theorem \ref{regularity}. 
Let 
\begin{equation*}
\begin{picture}(400,400)
\setsqparms[1`1`1`1;350`350]
\putVtriangle(0,00)%
[X` Z ` Y;g `f `h]
\end{picture}
\end{equation*}  
be the relative canonical model.  
Then the relative canonical K\"{a}hler current $\omega_{Z/Y}$ satisfies 
a partial differential equation on each fiber.  
Hence the regularity of $\omega_{Z/Y}$ (hence also the regularity of $d\mu_{can,X/Y}$) may be deduced by the parameter dependence of the solution of Monge-Amp\`{e}re equations. 

But after some time, I realized that this approach is extremely difficult 
to implement.
The reason is as follows.  Usually  since the canonical K\"{a}hler current 
is unique on each fiber of $h : Z\to Y$, it is natural to consider the 
variation of the canonical K\"{a}hler current satisfies a partial differential 
equation on each fiber which is (as is easily seen) 
essentially the Laplace equation with respect to the cacnonical K\"{a}hler
 current.  So far there is no difficulty.  The next step is to apply the 
 implicit function theorem.   Here the major difficulty arises.  
Because although the canonical K\"{a}hler current is $C^{\infty}$ on  
a nonempty Zariski open subset of 
each fiber, it is singular on a proper analytic subset of the fibers.
Hence it seems to be extremely difficult to fix the appropriate function space
 to apply the implicit function theorem.    Also it seems to be very difficult 
 to know the precise asymptotic behavior of the canonical K\"{a}hler current 
 near the singularities.  

Hence I decided to use the dynamical construction of canonical K\"{a}hler currents to deduce the (generic) horizontal smoothness of the relative canonical 
K\"{a}hler current.  

The advantage of this approach is that we can deduce the smoothness 
in terms of H\"{o}rmander's $L^{2}$-estimate for $\bar{\partial}$-operators. 
Because in each step, we only need to consider the variation of Bergman projections which is essentially a linear problem.  
In this way, we can deduce the regularity of the relative canonical K\"{a}hler current by the inductive estimates of Bergman projections.  

This inductive estimate is very smilar to the construction of Kuranishi family.

\subsection{Dynamical construction of the canonical K\"{a}hler currents}\label{Ds}

The canonical K\"{a}hler current  in Theorem \ref{main} 
can be constructed as the limit of a dynamical system as in (\cite{KE}).  

Let $X$ be a smooth projective $n$-fold with $\mbox{Kod}(X) \geqq 0$.
And let 
\begin{equation}
f : X -\cdots\rightarrow Y 
\end{equation}
be the Iitaka fibration associated with the complete linear system 
$|m_{0}!K_{X}|$ for some sufficiently large positive integer $m_{0}$.  By taking a suitable modifications, we  assume 
the followings:
\begin{enumerate}
\item[(1)] $Y$ is smooth and $f$ is a morphism.  
\item[(2)] $f_{*}\mathcal{O}_{X}(m_{0}!K_{X/Y})^{**}$ is a line bundle on $Y$, where
$**$ denotes the double dual.  
\end{enumerate}
We set 
\begin{equation}
L_{X/Y} := \frac{1}{m_{0}!}f_{*}\mathcal{O}_{X}(m_{0}!K_{X/Y})^{**}. 
\end{equation}  
In \cite{f-m} this $L_{X/Y}$ is denoted by $L_{X/Y}$. 
Let $a$ be  positive integer such that $f_{*}\mathcal{O}_{X}(aK_{X/Y})
\neq 0$. Then  we see that 
\begin{equation}\label{iso}
H^{0}(X,\mathcal{O}_{X}(maK_{X}))\simeq H^{0}(Y,\mathcal{O}_{Y}(ma(K_{Y}+L_{X/Y})))
\end{equation}
holds for every $m \geqq 0$.   In particular $\mbox{Kod}(X) = \dim Y$ holds. 
Hence by  (\ref{iso}), we see that $K_{Y} + L_{X/Y}$ is big. 
Let $A$ be an ample line bundle on $Y$ such that   
for every pseudoeffective singular hermitian line bundle $(F,h_{F})$ on 
$Y$, $\mathcal{O}_{Y}(jK_{Y}+ A + F)\otimes {\cal I}(h_{F})$ is globally generated for every $0\leqq j\leqq a$.

The existence of such an ample line bundle $A$ follows from  Nadel's vanishing theorem
(\cite[p.561]{n}).
Let $h_{A}$ be a $C^{\infty}$ hermitian metric on $A$ with strictly positive 
curvature. 
We  construct a sequence of singular hermitian metrics
$\{ h_{m}\}_{m\geqq 1}$ and a sequence of Bergman kernels $\{ K_{m}\}_{m\geqq 1}$ 
as follows. 

We set  
\begin{equation}
K_{1} := \left\{\begin{array}{ll} K(Y,K_{Y} +A,h_{A}), & \mbox{if}\,\, a > 1 \\ 
& \\ 
& \\
K(Y,K_{Y}+L_{X/Y}+A,h_{L_{X/Y}}\cdot h_{A}), & \mbox{if}\,\, a = 1 
\end{array}\right. 
\end{equation}
where for a singular hermitian line bundle $(F,h_{F})$ 
$K(Y,K_{Y}+F,h_{F})$ is (the diagonal part of) the Bergman kernel
of $H^{0}(Y,\mathcal{O}_{Y}(K_{Y} + F)\otimes {\cal I}(h_{F}))$ 
as (\ref{BergmanK}).
Then we set 
\begin{equation}
h_{1} := (K_{1})^{-1}. 
\end{equation}
We continue this process. 
Suppose that we have constructed 
$K_{m-1}$ and the singular hermitian metric $h_{m-1}$ on 
$(m-1)K_{Y}+\lfloor \frac{m-1}{a}\rfloor aL_{X/Y} + A$.  Then we define
\begin{equation}
K_{m}:= \left\{\begin{array}{ll} K(Y,mK_{Y}+\lfloor \frac{m}{a}\rfloor aL_{X/Y}+A,h_{m-1}) & \mbox{if}\,\, m  \not{\equiv} 0 \,\,\mbox{mod}\,a  \\
 & \\ 
 &  \\  
K(Y,m(K_{Y}+L_{X/Y})+A,h_{L_{X/Y}}^{a}\otimes h_{m-1}) & \mbox{if}\,\,m \equiv 0 \,\,\mbox{mod}\,a \end{array}\right.   
\end{equation}
and 
\begin{equation}
h_{m}:= (K_{m})^{-1}.  
\end{equation}
Thus inductively we construct the sequences $\{ h_{m}\}_{m\geqq 1}$
and $\{ K_{m}\}_{m \geqq 1}$.
This inductive construction is essentially the same one originated by the author in \cite{tu3}. 
The following theorem asserts that the above dynamical system yields the 
canonical K\"{a}hler current on $Y$. 

\begin{theorem}\label{DS}(\cite{canonical})
 Let $X$ be a smooth projective variety of nonnegative Kodaira dimension 
 and let $f : X \longrightarrow Y$ be the Iitaka fibration as above.  
 Let $m_{0}$ and $\{ h_{m}\}_{m \geq 1}$ be the sequence 
of hermitian metrics  as above and let $n$ denote $\dim Y$. 
Then 
\begin{equation}
h_{\infty} := \liminf_{m\rightarrow\infty} \sqrt[m]{(m!)^{n}\cdot h_{m}}
\end{equation}
is a singular hermitian metric on $K_{Y}+L_{X/Y}$ such that 
\begin{equation}
\omega_{Y}= \sqrt{-1}\,\Theta_{h_{\infty}}
\end{equation}
holds, where $\omega_{Y}$ is the canonical K\"{a}hler current on $Y$ 
as in Theorem \ref{main} and $n = \dim Y$.

More precisely 
\[
K_{\infty}:= \lim_{m\to\infty}h_{A}^{\frac{1}{am}}\cdot K_{am}^{\frac{1}{am}} 
\]
exists in $L^{1}$-topology (as a limit of bounded volume forms on Y) and 
$h_{\infty} = K_{\infty}^{-1}$ holds.  
In particular $\omega_{Y}= \sqrt{-1}\Theta_{h_{\infty}}$ (in fact $h_{\infty}$) is unique and is independent of the choice of $A$ and $h_{A}$. 
\fbox{}
\end{theorem} 
\begin{remark}
Similar theorem holds for a KLT pair with nonnegative (log) Kodaira dimension.
See \cite{LC}. But the corresponding dynamical system is not a 
single dynamical system, but is an infinite sequence of dynamical systems.  \fbox{}
\end{remark}
\subsection{Family of dynamical systems}
In this subsection, we shall consider the dynamical systems in Section \ref{Ds} on the relative canonical models. 
Let 
\begin{equation*}
\begin{picture}(400,400)
\setsqparms[1`1`1`1;350`350]
\putVtriangle(0,00)%
[X` Z ` Y;g `f `h]
\end{picture}
\end{equation*}  
be the relative canonical model as in Theorem \ref{regularity}.

Now we consider the relative version of the construction in Section \ref{Ds}.
Let $A$ be a sufficiently ample line bundle on $Z$ and let $h_{A}$ be 
a $C^{\infty}$-hermitian metric on $A$.
We set 
\begin{equation}
Y^{\circ}:= \{ y\in Y| \mbox{$f : X \to Y$ is smooth over $y$}\}.
\end{equation}
For every $y\in Y^{\circ}$ we  construct a sequence of singular hermitian metrics
$\{ h_{m,y}\}_{m\geqq 1}$ and a sequence of Bergman kernels $\{ K_{m,y}\}$ 
as follows. 

We set  
\begin{equation}
K_{1,y} := \left\{\begin{array}{ll} K(Z_{y},K_{Z_{y}} +A|Z_{y},h_{A}|Z_{y}), & \mbox{if}\,\, a > 1 \\ 
& \\ 
& \\
K(Z_{y},K_{Z_{y}}+L_{X/Z}|Z_{y}+A|Z_{y},h_{L_{X/Z}}\cdot h_{A}|Z_{y}), & \mbox{if}\,\, a = 1 
\end{array}\right. 
\end{equation}
Then we set 
\begin{equation}
h_{1,y} := (K_{1,y})^{-1}. 
\end{equation}
We continue this process. 
Suppose that we have constructed 
$K_{m-1}$ and the singular hermitian metric $h_{m-1}$ on 
$(m-1)K_{Z}+\lfloor \frac{m-1}{a}\rfloor aL_{X/Z} + A$.  Then we define
\begin{equation}
K_{m,y}:= \left\{\begin{array}{ll} K(Z_{y},mK_{Z_{y}}+\lfloor \frac{m}{a}\rfloor aL_{X/Z}|Z_{y}+A,h_{m-1,y}) & \mbox{if}\,\, m  \not{\equiv} 0 \,\,\mbox{mod}\,a \\
 & \\ 
 & \\  
K(Y,m(K_{Z_{y}}+L_{X/Z}|Z_{y})+A,h_{L_{X/Z}}^{a}|Z_{y}\otimes h_{m-1,y}) & \mbox{if}\,\,m \equiv 0 \,\,\mbox{mod}\,a \end{array}\right.  
\end{equation}
and 
\begin{equation}
h_{m,y}:= (K_{m,y})^{-1}.  
\end{equation}
Thus inductively we construct the sequences $\{ h_{m,y}\}_{m\geqq 1}$
and $\{ K_{m,y}\}_{m \geqq 1}$ for every $y\in Y^{\circ}$.

By \cite{b-p}, we see that 
$\sqrt{-1}\partial\bar{\partial}\log K_{m}$ extends to a closed positive 
current on $Y$.  Hence by Theorem \ref{DS}, we see that the relative canonical 
K\"{a}hler current: 
\begin{equation}
\omega_{Z/Y}:= \lim_{m\to\infty}\frac{\sqrt{-1}}{m}\partial\bar{\partial}\log K_{m} 
\end{equation}
extends to a closed positive current on $Y$.  
We denote the extended current again  by $\omega_{Z/Y}$. 
We shall prove Theorem \ref{regularity} by estimating the  
variation of $K_{m,y}$ with respect to the parameter $y\in Y^{\circ}$. 

%
%
%
%
  
\subsection{Variation of Bergman projections}

Let $U$ be an open subset of $Y^{\circ}$ such that 
$U$ is biholomorphic to the unit  polydisk $\Delta^{k}$ 
in $\mathbb{C}^{k}$ with ceneter $O$ via a local coordinate
 $(y_{1},\cdots,y_{k})$.
  Let $Z_{U} := h^{-1}(U)$ and let 
\begin{equation}
h_{U} : Z_{U} \to U
\end{equation}
be the restriction of $h$. 
Let us trivialize $h_{U} : Z_{U}\to U$ differentiably as 
\[
\Phi : Z_{U} \to Z_{0}\times U,
\]
where $Z_{0}$ denotes the central fiber $h_{U}^{-1}(O)$. 
Let 
\begin{equation}
P_{m,y} : L^{2}(Z_{y},A_{y} + mK_{Z_{y}}+\lfloor\frac{m}{a}\rfloor aL_{X/Z,y})
\to H^{0}(Z_{y},\mathcal{O}_{Z_{y}}(A_{y} + mK_{Z_{y}}+\lfloor\frac{m}{a}\rfloor aL_{X/Z,y}))
\end{equation}
the Bergman projection, i.e., the orthogonal projection   
with respect to the $L^{2}$-inner product $g_{m}$ defined by 
\begin{equation}\label{inner1}
g_{m}(\sigma_{y},\sigma^{\prime}_{y}) := 
\int_{Z_{y}}\sigma\cdot\overline{\sigma^{\prime}}\cdot h_{m-1,y}, 
\end{equation}
if $a\not{|}m$ and 
\begin{equation}\label{inner2}
g_{m}(\sigma_{y},\sigma^{\prime}_{y}) := 
\int_{Z_{y}}\sigma\cdot\overline{\sigma^{\prime}}\cdot 
h_{m-1,y}\cdot h_{X/Z,y}^{a},  
\end{equation}
if $a|m$.  Hearafter we shall omit $g_{m}$, if without fear of confusion.  
Then the above trivialization gives a trivialization:

\begin{equation}\label{2trivial}
\mathcal{L}^{2}(Z_{U},A|U + mK_{Z/Y}|U+\lfloor\frac{m}{a}\rfloor aL_{X/Z}|U)
\to L^{2}(Z_{0},A_{0} + mK_{Z_{0}}+\lfloor\frac{m}{a}\rfloor aL_{X/Z,0})
\times U,
\end{equation}
where $\mathcal{L}^{2}(Z_{U},A|U + mK_{Z/Y}|U+\lfloor\frac{m}{a}\rfloor aL_{X/Z}|U)$ denotes the Hilbert space bundle 

\begin{equation}\label{hilbbundle}
\pi_{L^{2}}: \mathcal{L}^{2}(Z_{U},A|U + mK_{Z/Y}|U+\lfloor\frac{m}{a}\rfloor aL_{X/Z}|U) \to U 
\end{equation}
such that 
\[
\pi_{L^{2}}^{-1}(y) := L^{2}(Z_{y},A_{y} + mK_{Z_{y}}+\lfloor\frac{m}{a}\rfloor aL_{X/Z}|Z_{y}). 
\]
Let 
\[
\bar{\partial}_{y} : C^{\infty}(Z_{y},A_{y} + mK_{Z_{y}}+\lfloor\frac{m}{a}\rfloor aL_{X/Z}|Z_{y}) \to 
A^{0,1}(Z_{y},A_{y} + mK_{Z_{y}}+\lfloor\frac{m}{a}\rfloor aL_{X/Z}|Z_{y})
\]
denote the $\bar{\partial}$-operator.
We set 
\begin{equation}
\mathbb{H}_{m,y} := H^{0}(Z_{y},\mathcal{O}_{Z_{y}}(A_{y} + mK_{Z_{y}}+\lfloor\frac{m}{a}\rfloor aL_{X/Z,y})). 
\end{equation}
Let $\sigma_{y}\in C^{\infty}(Z_{y},A_{y} + mK_{Z_{y}}+\lfloor\frac{m}{a}\rfloor aL_{X/Z}|Z_{y})$ be an arbitrary element. 
Let us consider the $\bar{\partial}$-equation:
\begin{eqnarray}
\bar{\partial}_{y}(Q_{m,y}(\sigma_{y}))  & =  & \bar{\partial}_{y}\sigma_{y} \\
Q_{m,y}(\sigma_{y}) & \perp & \mathbb{H}_{m,y}. 
\end{eqnarray}
Then 
\begin{equation}
Q_{m,y} : \mathbb{L}^{2}_{m} \to \mathbb{H}_{m,y}^{\perp} 
\end{equation}
is the orthogonal projection. 
Hence the Bergman projection is given by 
\begin{equation}
P_{m,y}(\sigma_{y}) = \sigma_{y} - Q_{m,y}(\sigma_{y}). 
\end{equation}
This implies that 
\begin{equation}
D_{y}P_{m,y} = -D_{y}Q_{m,y}
\end{equation}
holds, where $D_{y}$ denotes the hermitian connection with respect to 
$g_{m}$ (cf. (\ref{inner2})). 

Let us calculate the variation of $P_{m,y}$ at $y = 0$. 
Let $\sigma \in \mathbb{H}_{m,0}$ and let us extend $\sigma$ 
as a section $\tilde{\sigma}$ of the Hilbert space bundle (\ref{hilbbundle}) 
by the parallel displacement with respect to the hermitian connection with respect to $g_{m}$ along a  smooth curve on $Y$.  We note that since the connection may not be flat, the parallel displacement depends on the choice of the smooth 
curve.  Hereafter we shall fix a differential curve to fix the extension
 $\tilde{\sigma}$. 

Then differentiating the equation: 
\[
\bar{\partial}_{y}\tilde{\sigma} (y) = \bar{\partial}_{y}Q_{m,y}(\tilde{\sigma} (y))
\]
with respect to $y$ at $y = 0$, we obtain the equation:
\begin{equation}\label{dbar}
\theta_{m,0}(\sigma) = \bar{\partial}_{0}(D_{y}Q_{m,y}(\sigma))
\end{equation}
where $\theta_{m,0}$ represents the Kodaira-Spencer class.

We shall decompose $D_{y}Q_{m,y}(\tilde{\sigma})$ as 
\begin{equation}
D_{y}Q_{m,y}(\tilde{\sigma}) = D_{y}Q_{m,y}(\tilde{\sigma})_{\mathbb{H}}
 + D_{y}Q_{m,y}(\tilde{\sigma})_{\mathbb{H}^{\perp}} 
\end{equation}
corresponding to the orthogonal decomposition 
\[
\mathbb{L}_{m,y} = \mathbb{H}_{m,y} \oplus \mathbb{H}_{m,y}^{\perp}. 
\]
Then we have that 
\begin{equation}\label{dbar2}
\theta_{m,0}(\sigma) = \bar{\partial}_{0}(D_{y}Q_{m,y}(\sigma))_{\mathbb{H}^{\perp}}
\end{equation}
holds, i.e., $D_{y}Q_{m,y}(\sigma)_{\mathbb{H}^{\perp}}$ is the minimal 
solution of (\ref{dbar2}). 
Now we shall fix the standard K\"{a}hler metric on $Y \sim \Delta$ 
induced by the standard K\"{a}hler metric on $\mathbb{C}$. 
Let us estimate the operator norm of  
\begin{equation}
(D_{y}Q_{m,y})_{\mathbb{H}^{\perp}} : \mathbb{H}_{m,0}\to \mathbb{H}_{m,0}^{\perp}. 
\end{equation}
The norm is estimated by H\"{o}rmander's $L^{2}$-estimate for  
$\bar{\partial}$-operators.

First we see that $\theta_{m,0}$ consists of the Kodaira-Spencer class 
of the deformation of $Z_{m,0}$ and the Kodaira-Spencer class of 
the bundle $\lfloor m/a\rfloor aL_{X/Z,0}$. 
Then there exists a positive constant $C_{0}$ independent of $m$ such that 
\begin{equation}
\parallel\theta_{m,0}\parallel_{L^{\infty}} \leqq C_{0}
\end{equation} 
holds. 
On the other hand by H\"{o}rmander's $L^{2}$-estimate, we see that 
\begin{lemma}\label{norm} There exists a positive constant $C_{1}$ such that 
\begin{equation}
\parallel (D_{y}Q_{m,y})_{\mathbb{H}^{\perp}}\parallel \leqq C_{1}
\end{equation}
holds for every $m \geqq 1$. \fbox{}
\end{lemma}
For $k=2$, differentiating (\ref{dbar2}), we have the equation
\begin{equation}
(D_{y}\theta_{m,y})(\sigma)(0) = \bar{\partial}_{0}(D_{y}^{2}Q_{m,y}(\sigma)_{\mathbb{H}^{\perp}})+ \theta_{m,0}(D_{y}Q_{m,y}(\sigma)_{\mathbb{H}^{\perp}}). 
\end{equation}
Hence we may estimate $D_{y}^{2}Q_{m,y}(\sigma)_{\mathbb{H}^{\perp}}$ 
as 
\[
\parallel D_{y}^{2}Q_{m,y}(\sigma)_{\mathbb{H}^{\perp}}\parallel
\leqq C_{2}
\]
for some positive constant $C_{2}$ independent of $m$. 
For $k \geqq 2$, inductively we have:
\begin{lemma}
For every $k\geqq 1$, there exists a positive constant $C_{k}$ such that 
\begin{equation}
\parallel (D_{y}^{k}Q_{m,y})_{\mathbb{H}^{\perp}}\parallel \leqq C_{k}
\end{equation}
holds for every $m \geqq 1$. \fbox{}
\end{lemma}

\subsection{Estimate of the holomorphic part}
Now we shall estimate the holmorphic part of the derivatives of $Q_{m,y}$ at $y= 0$.
Let 
\[
\tau_{hol} \in H^{0}(Z,\mathcal{O}_{Z}(mK_{Z/Y}+\lfloor \frac{m}{a}\rfloor aL_{X/Z})),
\] 
be an arbitrary holomorphic section. 
Differentiating the trivial identity:
\begin{equation}
g_{m}(Q_{m,y}(\tilde{\sigma}),\tau_{hol}) = 0,
\end{equation}  
we obtain that for every positive integer $\ell$ 
\begin{equation}\label{induc}
\sum_{i+j= \ell}\int_{X}h_{m-1}\left(D_{y}^{i}Q_{m,y}(\tilde{\sigma}),D_{y}^{j}\tau_{hol}\right) = 0
\end{equation}
holds, where $g_{m}$ denotes the $L^{2}$-metric defined by (\ref{inner1}) and 
(\ref{inner2}). 
Hence (\ref{induc}) implies that we can estimate 
$D^{k}_{y}Q_{m,y}$ in terms of the esimate of 
$\{ D_{y}^{\ell}Q_{m-1,y}\}_{\ell = 0}^{k-1}$.  

\begin{lemma}\label{norm2} There exists a positive constant $C^{\prime}_{k}$ 
independent of $m$ such that 
\begin{equation}
\parallel (D^{k}_{y}Q_{m,y})_{\mathbb{H}}\parallel  \leqq C_{k}^{\prime}
\end{equation}
holds for every $m \geqq 1$. \fbox{}
\end{lemma}
{\em Proof}.
Let $\tau_{0}$ be an element of
 $H^{0}(Z_{0},\mathcal{O}_{Z_{0}}(mK_{Z/Y}+\lfloor \frac{m}{a}\rfloor aL_{X/Z}))$.
We extend $\tau_{0}$ to the $\ell$-th infinitesimal neighbourhood $Z_{0}^{(\ell)}$ of $Z_{0}$ by the successive extension.    
By the $L^{2}$-estimates, we may take the extension $\tau^{(\ell)}$ so that 
\begin{equation}
\parallel D^{\ell}\tau^{(\ell)}\parallel_{(\ell)} \leqq C_{(\ell)}
\end{equation} 
holds for some positive constant $C_{(\ell)}$ independent of $m$.
Then replacing $\tau_{hol}$ by $\tau^{(\ell)}$ in  (\ref{induc}), by induction on $\ell$, we have the lemma. \fbox{} 

\subsection{Variation of Bergman kernels}

The parameter dependence of Bergman kernels can be deduced from 
the variation of the Bergman projections.
By the trivial equality:
\begin{equation}
(D_{y}^{k}P_{m})(\tilde{\sigma})(z)
= \int_{X(\zeta)} h_{m-1}\cdot D_{y}^{k}K_{m}(z,\zeta )\cdot\tilde{\sigma}(\zeta),
\end{equation}
(where the integral is taken with respect to the parameter $\zeta$)
\begin{equation}
\sum_{i=0}^{k}(D_{y}^{k-i}K_{m}(z,w),D_{y}^{i}Q_{m}(\tilde{\sigma}))_{g_{m}} = 0
\end{equation}
holds. 
Then by induction on $k$ and the extremal property of Bergman kernels, there exists a 
positive constant $C$ independent of $m$ such that 
\begin{equation}\label{single}
|(D_{y}^{k}P_{m})(\tilde{\sigma})(z)|_{h_{m-1}}
= \left(D_{y}^{k}K_{m}(z,\zeta ),\tilde{\sigma}(\zeta)\right)\leqq C\cdot m^{\frac{n}{2}}
\parallel\tilde{\sigma}\parallel
\end{equation}
holds for every $z\in Z_{y}$. 
Combining (\ref{single}), this implies that there exists a positive constant $C(k)$ depending only 
on $k$ such that 
\begin{equation}
 |m^{-n}D_{y}^{k}K_{m}(z,\zeta)|_{h_{m-1}} < C(k)
\end{equation}
holds on $Z_{y}$. 
Hence by the Sobolev's embedding theorem, we see that 
there exists a positive constant $\hat{C}_{k}$ independent of $m$ such that 
\begin{equation}
|\left((m!)^{-n}K_{m}\right)^{\frac{1}{m}}|_{C^{k}} \leqq \hat{C}_{k}
\end{equation}
holds on $Z_{y}$ 
By Theorem \ref{DS}, this means that the relative canonical measure $d\mu_{can,X/Y}$ is $C^{\infty}$ on a nonempty Zariski open subset of $X$. 
This completes the proof of Theorem \ref{regularity}. \fbox{} 



\section{Appendix}

In this section, we collect several analytic tools used in this article. 

\subsection{Ampleness criterion for  line bundles on quasiprojective varieties}
In this section we prove a criterion of quasiprojectivity 
used in the previous section.   The criterion is almost the same as 
in \cite{sch-t}.  But it is slightly stronger.  

Let $X$ be a not necessarily reduced algebraic space with
compactification $\ol{X}$ in the sense of algebraic spaces, and
let $L$ be a holomorphic line bundle on $\ol{X}$ with a {\it
positive} singular hermitian metric $h$ in the following sense. 

\begin{definition}\label{singdef} Let $Z$ be a reduced complex space
and $L$ a holomorphic line bundle. A {\em singular} hermitian
metric $h$ on $L$ is a singular hermitian metric $h$ on
$L|Z_{reg}$ with the following property: There exists a
desingularization $\pi: \wt{Z} \to Z$ such that $h$ can be
extended from $Z_{reg}$ to a singular hermitian metric $\wt{h}$
on $\pi^*L$ over $\wt{Z}$. \fbox{}
\end{definition}

\begin{condition}[\bf P]
We say that the positivity condition (P) holds, if
\begin{enumerate}
\item[(i)] For all $p\in  X$ and any holomorphic curve $C \subset
\ol X$ through $p$ the (positive,
$d$-closed) current $\sqrt{-1}\Theta_h|C$ is well-defined, and the Lelong
number $\nu(\sqrt{-1}\Theta_h|C,p)$ vanishes,
\item[(ii)] For any smooth locally closed subspace $Z \subset X$
of $\dim Z > 0$, $h|_{Z}$ is well defined and 
$(L|_{Z},h|_{Z})$ is big (cf. Definition \ref{numerical} below).  \fbox{}
\end{enumerate}
\end{condition}

Now we state the criterion.

\begin{theorem}\label{criterion}
Let $X$ be an irreducible, not necessarily reduced algebraic space
with a compactification $\ol X$. Let $L$ be a holomorphic line
bundle on $\ol X$. The map
$$
\Phi_{|m L|} : \ol X \rightharpoonup \mathbb P^N(m),
$$
where $N(m)=\dim |mL|$, defines an embedding of $X$ for
sufficiently large $m$, if it satisfies condition (P). \fbox{}
\end{theorem}

The proof of Theorem \ref{criterion} is essentially the same as 
the one of \cite[Theorem 6]{sch-t} except the use of Theorem \ref{kodaira} below to perturbe the metric to a metric with strictly positive curvature.
 
\subsection{Kodaira's lemma for big pseudoeffective line bundles}

In this subsection, we prove a singular hermitian version of Kodaira's lemma
(cf. \cite[Appendix]{k-o}). 

First we shall define the big singular hermitian line bundle.  

\begin{definition}\label{numerical}
$(L,h_{L})$ be a pseudoeffective singular hermitian line bundle on a projective  manifold  $X$.  
We set 
\begin{equation}
\nu_{num}(L,h_{L}) := \sup \{\dim V \mid \mbox{\em $V$ is a subvariety of $X$ such that 
$h_{L}\!\mid_{V}$ is well defined} 
\end{equation}
\begin{equation}
\hspace{30mm}\mbox{\em and $(L,h_{L})^{\dim V}\!\!\cdot V> 0$}\}. 
\end{equation}  
We call $\nu_{num}(L,h_{L})$ the {\bf numerical Kodaira dimension} of $(L,h_{L})$.  If $\nu_{num}(L,h_{L}) = \dim X$ we say that $(L,h_{L})$ is big. 
  \fbox{}
\end{definition}

\begin{lemma}\label{generic}
Let $X$ be a smooth projective variety and let $|H|$ be a very ample linear 
system.  
Then there  exists a smooth member $H^{\prime}\in \mid\! H\!\mid$, such that 
\begin{equation}
{\cal I}(h_{L}^{m})\otimes\mathcal{O}_{H^{\prime}}
= {\cal I}(h_{L}^{m}\mid_{H^{\prime}})
\end{equation}
holds for every $m\geqq 1$. \fbox{}   
\end{lemma}
{\em Proof of Lemma \ref{generic}}.
Let $A$ be a sufficiently ample line bundle such that 
$\mathcal{O}_{X}(A+mL))\otimes {\cal I}(h_{L}^{m})$  is globally generated
for all $m\geqq 1$.  
Let $\{ \sigma^{(m)}_{j}\}_{j=1}^{N_{m}}$ be a (complete) basis 
of $H^{0}(X,\mathcal{O}_{X}(A+mL))\otimes {\cal I}(h_{L}^{m}))$. 
We consider the subset
\begin{equation}
U:= \{ F\in \mid\! H\!\mid ; \mbox{$F$ is smooth}, 
\int_{F}\mid\sigma_{j}^{(m)}\mid^{2}\cdot h_{L}^{m}\cdot h_{A}\cdot \, dV_{F}
< + \infty \,\,\,
\end{equation}
\vspace{-8mm}
\begin{equation}
\hspace{50mm}\mbox{for every $m$ and $1\leqq j\leqq N_{m}$.}\}
\end{equation}
of $\mid\! H\!\mid$, where $dV_{F}$ denotes the volume form  on $F$ induced by the 
K\"{a}hler form $\omega$. 
We claim that such $U$ is the complement of at most a countable union of 
proper subvarieties of $\mid\! H\!\mid$.
Let us fix a positive integer $m$.  
\begin{equation}
E_{m}:= \{ F\in \mid\! H\!\mid ; \mbox{$F$ is smooth}, 
\int_{F}\mid\sigma_{j}^{(m)}\mid^{2}\cdot h_{L}^{m}\cdot h_{A}\,\, dV_{F}
=  + \infty \}
\end{equation}
is of measure $0$ by Fubini's theorem.  
Then since $U = |H| - \cup_{m=1}^{\infty}E_{m}$, we complete the proof of 
Lemma \ref{generic}.   
\fbox{} 
 
\begin{theorem}\label{kodaira}
Let $X$ be a projective manifold and let $(L,h_{L})$ be a big 
psedoeffective singular hermitian line bundle. 
Then there exists a singular hermitian metric $h^{+}_{L}$ on $L$ such that 
\begin{enumerate}
\item[(1)] $\sqrt{-1}\Theta_{h^{+}_{L}}$ is strictly positive 
everywhere on $X$, 
\item[(2)] $h^{+}_{L} \geqq h_{L}$ holds on $X$. 
\end{enumerate} \fbox{}
\end{theorem}
\noindent 
Let us explain the relation between Theorem\ref{kodaira} and 
the original Kodaira's lemma .
Let $D$ be an ample divisor on a smooth projective variety $X$.   Let us identify divisors with line bundles. 
By Kodaira's lemma,  there exists a $C^{\infty}$ hermitian metrics 
$h_{D},h_{E}$ on $D,E$ respectively (the notion of hermitian metrics 
naturally extends to the case of $\mathbb{Q}$-line bundles) such that 
the curvature of $h_{D}\cdot h_{E}^{-1}$ is stricly positive. 
Let $\sigma_{E}$ be a multivalued holomorphic section of $E$ with divisor $E$
such that $h_{E}(\sigma_{E},\sigma_{E})\leqq 1$ on $X$. 
Then 
\begin{equation}
h_{D}^{+}:= \frac{h_{D}}{h_{E}(\sigma_{E},\sigma_{E})}
\end{equation}
is a singular hermitian metric on $D$  such that 
\begin{enumerate}
\item[(1)] $\sqrt{-1}\Theta_{h_{D}^{+}}$
is strictly positive everywhere on $X$.
\item[(2)] $h_{D} \leqq h_{D}^{+}$
holds on $X$.
\end{enumerate}
In this way Theorem \ref{kodaira} can be viewed as an analogue of the usual 
Kodaira's lemma to the case of big pseudoeffective singular hermitian line bundles.
\subsection{Proof of Theorem \ref{kodaira}}
The proof of Theorem \ref{kodaira} presented here is not very much different from the original proof of Kodaira's lemma (cf. \cite{ka1} or \cite[Appendix]{k-o}).  But it requires estimates of  Bergman kernels  and  additional care for the multiplier ideal sheaves. 

Let $X$ be a smooth projective variety of dimension $n$ and let 
$(L,h_{L})$ be a big pseudoeffective singular hermitian line bundle on $X$. 
Let $\omega$ be a K\"{a}hler form on $X$ and let $dV$ be the associated 
volume form on $X$. 
Let $H$ be a smooth very ample divisor on $X$.  
\begin{lemma}\label{big}
There exists a positive integer $m_{0}$ such that 
$m_{0}(L,h_{L}) - H$ is big, i.e.,  
\begin{equation}
\limsup_{\ell\rightarrow\infty}\ell^{-n}\cdot\dim H^{0}(X,\mathcal{O}_{X}(\ell (m_{0}L - H)\otimes {\cal I}(h_{L}^{m_{0}\ell})) > 0
\end{equation}
holds.  \fbox{}
\end{lemma}
{\em Proof of Lemma \ref{big}.}
Replacing $H$ by a suitable member of $\mid\! H\!\mid$, by Lemma \ref{generic},
we may assume that 
\begin{equation}
{\cal I}(h_{L}^{m})\mid_{H} = {\cal I}(h_{L}^{m}\mid_{H})
\end{equation}
holds for every $m \geqq 1$. 
Let us consider the exact sequence 
\begin{equation}
0\rightarrow H^{0}(X,\mathcal{O}_{X}(mL - H)\otimes {\cal I}(h_{L}^{m}))\rightarrow H^{0}(X,\mathcal{O}_{X}(mL)\otimes {\cal I}(h_{L}^{m})) 
\end{equation}
\begin{equation}
\hspace{45mm}\rightarrow H^{0}(H,\mathcal{O}_{H}(mL)\otimes {\cal I}(h_{L}^{m}\mid_{H})). 
\end{equation}
Then since $\mu (L,h_{L}) > 0$ and 
\begin{equation}
\dim H^{0}(H,\mathcal{O}_{H}(mL)\otimes {\cal I}(h_{L}^{m}\mid_{H})) = O(m^{n-1})
\end{equation}
we see that for every sufficiently large $m$, 
\begin{equation}
H^{0}(X,\mathcal{O}_{X}(mL - H)\otimes {\cal I}(h_{L}^{m}))\neq 0
\end{equation}
holds. \fbox{}  \vspace{3mm}\\
To prove Lemma \ref{big}, we need to refine the above argument a little bit. 
Let $m_{0}$ be a positive integer such that 
\begin{equation}\label{m}
m_{0} >n\cdot \frac{(L,h_{L})^{n-1}\!\!\cdot H}{(L,h_{L})^{n}} 
\end{equation}
holds.   For very general $H_{1}^{(\ell)},\cdots H^{(\ell)}_{\ell}\in \mid\! H\!\mid$, by Lemma \ref{generic}, replacing $m$ by $m_{0}\ell$ and 
$H$ by $\ell H$, we have the exact sequence 
\begin{equation}
0\rightarrow H^{0}(X,\mathcal{O}_{X}(\ell (m_{0}L - H))\otimes{\cal I}(h_{L}^{m_{0}\ell}))\rightarrow H^{0}(X,\mathcal{O}_{X}(m_{0}\ell L)\otimes {\cal I}(h_{L}^{m_{0}\ell})).  
\end{equation}
\begin{equation}
\hspace{40mm} \rightarrow 
\oplus_{i=1}^{\ell}H^{0}(H^{(\ell)}_{i},\mathcal{O}_{H_{i}}(m_{0}\ell L)\otimes {\cal I}(h_{L}^{m_{0}\ell}\mid_{H_{i}})). 
\end{equation}
We note that $\{ H_{i}^{(\ell)}\}_{i=1}^{\ell}$ are  chosen for each $\ell$.
If we take $\{ H_{i}^{(\ell)}\}_{i=1}^{\ell}$ very general, 
we may assume  that 
\begin{equation}
\dim H^{0}(H^{(\ell)}_{i},\mathcal{O}_{H_{i}}(mL)\otimes {\cal I}(h_{L}^{m}\mid_{H_{i}}))
\end{equation}
is independent of $1\leqq i\leqq \ell$ for every $m$. 
This implies that 
\begin{equation}
\limsup_{\ell\rightarrow\infty}\ell^{-n}\cdot\dim H^{0}(X,\mathcal{O}_{X}(\ell (m_{0}L - H))\otimes {\cal I}(h_{L}^{m_{0}\ell}))
\end{equation}
\begin{equation}
\hspace{20mm}
\geqq \frac{1}{n!}(L,h_{L})^{n}\cdot m_{0}^{n}- \frac{1}{(n-1)!}
 \{(L,h_{L})^{n-1}\!\!\cdot H\} \cdot m_{0}^{n-1}
\end{equation}
holds.  By (\ref{m}), we see that 
\begin{equation}
\frac{1}{n!}(L,h_{L})^{n}\cdot m_{0}^{n}- \frac{1}{(n-1)!}\{(L,h_{L})^{n-1}\!\!\cdot H\}m_{0}^{n-1}
\end{equation}
is positive. 
This completes the proof of Lemma \ref{big}. \fbox{}
\vspace{5mm}\\
Let $A$ be a sufficiently ample line bundle on $X$ and let $h_{A}$ be a 
$C^{\infty}$ hermitian metric such that the curvature of $h_{A}$ is everywhere
strictly positive on $X$.  Here the meaning of ``sufficiently ample'' will 
be specified later. 
Let $m$ be a positive integer. 
Let us consider the inner product 
\begin{equation}
(\sigma ,\sigma^{\prime}) := \int_{X}h_{A}\cdot h_{L}^{m}\cdot \sigma
\cdot\bar{\sigma}^{\prime}\, dV
\end{equation}
on $H^{0}(X,\mathcal{O}_{X}(A+ mL)\otimes {\cal I}(h_{L}^{m}))$ and 
let $K_{m}$ be the associated (diagonal part of) Bergman kernel.
Let us consider the subspace: 
\begin{equation} 
H^{0}(X,\mathcal{O}_{X}(A+\ell (m_{0}L - H))\otimes {\cal I}(h_{L}^{m_{0}\ell}))
\subset H^{0}(X,\mathcal{O}_{X}(A+m_{0}\ell L)\otimes {\cal I}(h_{L}^{m_{0}\ell}))
\end{equation}
as a Hilbert subspace and let $K_{m_{0}\ell}^{+}$ denotes the associated 
Bergman kernel with respect to the restriction of the inner product 
on \\ $H^{0}(X,\mathcal{O}_{X}(A + m_{0}\ell L)\otimes {\cal I}(h_{L}^{m_{0}\ell}))$
to the subspace
$H^{0}(X,\mathcal{O}_{X}(A+\ell (m_{0}L - H))\otimes {\cal I}(h_{L}^{m_{0}\ell}))
$. 
Then by definition, we have the trivial inequality :
\begin{equation}\label{trivial inequality}
K_{m_{0}\ell}^{+} \leqq K_{m_{0}\ell}
\end{equation}
holds on $X$ for every $\ell \geqq 1$. 

The next lemma follows from the same argument as in \cite{dem}

\begin{lemma}\label{reciprocity}(\cite{dem})
If $A$ is sufficiently ample,
\begin{equation}
h_{L} := \mbox{the lower envelope of}\,\,(\limsup_{m\rightarrow\infty}\sqrt[m]{K_{m}})^{-1}.
\end{equation}
holds. \fbox{}
\end{lemma}
\begin{remark}
In \cite{dem}, Demailly considered the local version of Lemma \ref{reciprocity}, but the same proof works thanks to the sufficiently ample line bundle $A$. \fbox{}
\end{remark}
We note that 
\begin{equation}
\int_{X}h_{A}\cdot h_{L}^{m}\cdot K_{m}\cdot dV = \dim H^{0}(X,\mathcal{O}_{X}(A+ mL)\otimes {\cal I}(h_{L}^{m}))
\end{equation}
and 
\begin{equation}
\int_{X}h_{A}\cdot h_{L}^{m_{0}\ell}\cdot K_{m_{0}\ell}^{+}\cdot dV = \dim H^{0}(X,
\mathcal{O}_{X}(A+\ell (m_{0}L - H))\otimes {\cal I}(h_{L}^{m_{0}\ell}))
\end{equation}
hold.  
Hence by Lemma \ref{big}
\begin{equation}\label{ve}
\limsup_{\ell\rightarrow\infty}\, \left((m_{0}\ell)^{-n}\cdot \int_{X}h_{L}^{m_{0}\ell}\cdot K^{+}_{m_{0}\ell}\cdot dV\right)  > 0
\end{equation}
holds. 
Then by Fatou's lemma, we see that 
\begin{equation}
\int_{X}\limsup_{\ell\rightarrow\infty}\frac{h_{A}\cdot h_{L}^{m_{0}\ell}\cdot K^{+}_{m_{0}\ell}}{(m_{0}\ell )^{n}}
\geqq
\limsup_{\ell\rightarrow\infty}\int_{X}\frac{h_{A}\cdot h_{L}^{m_{0}\ell}\cdot K^{+}_{m_{0}\ell}}{(m_{0}\ell)^{n}} > 0 
\end{equation} 
hold.
In particular
\begin{equation}
\limsup_{\ell\rightarrow\infty}\frac{h_{A}\cdot h_{L}^{m_{0}\ell}\cdot K^{+}_{m_{0}\ell}}{(m_{0}\ell )^{n}}
\end{equation}
is not identically $0$\footnote{At this moment, there is a possibility that it is 
identically $+\infty$.}.
This implies that 
\begin{equation}
\limsup_{\ell\rightarrow\infty}\sqrt[m_{0}\ell]{K^{+}_{m_{0}\ell}}
\end{equation}
is not identically $0$ and by Lemma \ref{reciprocity} and (\ref{trivial}), 
it is finite.  
\noindent Let $h_{H}$ be a $C^{\infty}$ hermitian metric on $H$ with strictly positive
curvature and let $\tau$ be a global holomorphic section of $\mathcal{O}_{X}(H)$
with divisor $H$ such that $h_{H}(\tau ,\tau) \leqq 1$ holds on $X$. 
We set 
\begin{equation}
h_{L}^{+} := (\limsup_{\ell\rightarrow\infty}\sqrt[m_{0}\ell]{K^{+}_{m_{0}\ell}}\,\,)^{-1}\cdot h_{H}(\tau ,\tau ).
\end{equation}
Then $h_{L}^{+}$ is a singular hermitian metric
on $L$, since \\ $(\limsup_{\ell\rightarrow\infty}\sqrt[m_{0}\ell]{K^{+}_{m_{0}\ell}}\,\,)^{-1}\cdot \mid\tau\mid^{2}$ can be viewed as  a singular hermitian metric 
on $L - H$ with semipositive curvature current.  By the construction it is clear that the curvature current of $h_{L}^{+}$ is bigger than or equal to the curvature of $h_{H}$. 
In particular the curvature current of $h_{L}^{+}$ is strictly positive. 
And by the construction 
\begin{equation}
h_{L} \leqq h_{L}^{+}
\end{equation}
holds on $X$. 
This completes the proof of Theorem \ref{kodaira}. \fbox{}

Author's address\\
Hajime Tsuji\\
Department of Mathematics\\
Sophia University\\
7-1 Kioicho, Chiyoda-ku 102-8554\\
Japan \\
e-mail address: tsuji@mm.sophia.ac.jp  or h-tsuji@h03.itscom.net
\end{document}